\documentclass[11pt,a4paper]{scrartcl}

\usepackage[a4paper,verbose]{geometry}
\pdfoutput=1

\usepackage[utf8]{inputenc}
\usepackage[T1]{fontenc} %
\usepackage{libertine}

\usepackage{amsfonts}
\usepackage{amssymb}
\usepackage{amsthm}
\usepackage{amsmath}
\usepackage{caption}
\usepackage{etoolbox}
\usepackage{stmaryrd}
\usepackage{ifthen}
\usepackage{suffix}
\usepackage{verbatim} %

\usepackage{xfrac}

\usepackage{bm}
\usepackage{physics}
\usepackage[safe]{tipa} %
\usepackage{mathtools} %
\usepackage{color}
\usepackage{wrapfig}
\usepackage{tikz}
\usepackage{tikzit}
\usetikzlibrary{cd, babel}

\usepackage[framemethod=tikz]{mdframed}

\usepackage[backend=biber,style=numeric-comp,sorting=none,date=short,maxcitenames=2]{biblatex}

\usepackage{hyperref}

\newcommand*{\email}[1]{%
  \normalsize\href{mailto:#1}{\texttt{#1}}\par
}

\bibliography{bibliography}

\newcommand*{\citep}[1]{\parencite{#1}}
\newcommand*{\citet}[1]{\textcite{#1}}

\usepackage[capitalize]{cleveref} %

\usepackage{graphicx}
\usepackage{pict2e}

\def\mdash{{\hbox{-}}}

\makeatletter
\newcommand{\adjunction}{\@ifstar\named@adjunction\normal@adjunction}
\newcommand{\normal@adjunction}[4]{%
  #1\colon #2%
  \mathrel{\vcenter{%
    \offinterlineskip\m@th
    \ialign{%
      \hfil$##$\hfil\cr
      \longrightharpoonup\cr
      \noalign{\kern-.3ex}
      \smallbot\cr
      \longleftharpoondown\cr
    }%
  }}%
  #3 \noloc #4%
}
\newcommand{\named@adjunction}[4]{%
  #2%
  \mathrel{\vcenter{%
    \offinterlineskip\m@th
    \ialign{%
      \hfil$##$\hfil\cr
      \scriptstyle#1\cr
      \noalign{\kern.1ex}
      \longrightharpoonup\cr
      \noalign{\kern-.3ex}
      \smallbot\cr
      \longleftharpoondown\cr
      \scriptstyle#4\cr
    }%
  }}%
  #3%
}
\newcommand{\longrightharpoonup}{\relbar\joinrel\rightharpoonup}
\newcommand{\longleftharpoondown}{\leftharpoondown\joinrel\relbar}
\newcommand\noloc{%
  \nobreak
  \mspace{6mu plus 1mu}
  {:}
  \nonscript\mkern-\thinmuskip
  \mathpunct{}
  \mspace{2mu}
}
\newcommand{\smallbot}{%
  \begingroup\setlength\unitlength{.15em}%
  \begin{picture}(1,1)
  \roundcap
  \polyline(0,0)(1,0)
  \polyline(0.5,0)(0.5,1)
  \end{picture}%
  \endgroup
}
\makeatother

\let\op=\relax

\DeclareMathOperator{\op}{^\text{op}}

\newcommand{\nn}{{\mathbb{N}}}
\newcommand{\rr}{{\mathbb{R}}}
\newcommand{\Tt}{{\mathbb{T}}}

\newcommand{\Cat}[1]{\mathbf{#1}}
\newcommand{\cat}[1]{\mathcal{#1}}
\newcommand{\Fun}[1]{\mathsf{#1}}

\newcommand{\Kl}{\mathcal{K}\mspace{-2mu}\ell}

\DeclareMathOperator*{\E}{\mathbb{E}}

\renewcommand{\d}{\mathrm{d}}

\DeclareMathOperator{\Pa}{\mathcal{P}}

\DeclareMathOperator{\id}{\mathsf{id}}

\newcommand{\xto}[1]{\xrightarrow{#1}}

\makeatletter
\newcommand{\mathoverlap}[2]{\mathpalette\mathoverlap@{{#1}{#2}}}
\newcommand{\mathoverlap@}[2]{\mathoverlap@@{#1}#2}
\newcommand{\mathoverlap@@}[3]{\ooalign{$\m@th#1#2$\crcr\hidewidth$\m@th#1#3$\hidewidth}}
\makeatother

\newcommand{\klcirc}{\bullet} %

\makeatletter
\providecommand*{\xmapstofill@}{%
  \arrowfill@{\mapstochar\relbar}\relbar\rightarrow
}
\providecommand*{\xmapsto}[2][]{%
  \ext@arrow 0395\xmapstofill@{#1}{#2}%
}

\makeatletter
\def\slashedarrowfill@#1#2#3#4#5{%
  $\m@th\thickmuskip0mu\medmuskip\thickmuskip\thinmuskip\thickmuskip
   \relax#5#1\mkern-7mu%
   \cleaders\hbox{$#5\mkern-2mu#2\mkern-2mu$}\hfill
   \mathclap{#3}\mathclap{#2}%
   \cleaders\hbox{$#5\mkern-2mu#2\mkern-2mu$}\hfill
   \mkern-7mu#4$%
}
\def\rightslashedarrowfill@{%
  \slashedarrowfill@\relbar\relbar\mapstochar\rightarrow}
\newcommand\xslashedrightarrow[2][]{%
  \ext@arrow 0055{\rightslashedarrowfill@}{#1}{#2}}
\makeatother

\theoremstyle{definition}
\newtheorem{defn}{Definition}[section]

\newtheorem{ex}[defn]{Example}

\newtheorem{rmk}[defn]{Remark}

\newtheorem{question}[defn]{Question}

\newtheorem{prop}[defn]{Proposition}
\newtheorem{prop*}{Proposition}

\newtheorem{lemma}[defn]{Lemma}

\newtheorem{cor}[defn]{Corollary}

\newtheorem*{thm*}{Theorem}
\newtheorem*{cor*}{Corollary}

\definecolor{darkblue}{rgb}{0,0,0.7} 
\newcommand{\red}[1]{{\color{red} #1}}

\usepackage{authblk}

\author{Toby St. Clere Smithe
  \thanks{Early draft of work in progress. Comments and suggestions welcome.}}
\affil{Topos Institute \\ \email{toby@topos.institute}}

\hypersetup{
  pdftitle={Dynamics over Polynomials},
  pdfauthor={Toby St Clere Smithe}
}

\tikzstyle{xshiftu}=[shift = {(#1, 0)}]
\tikzstyle{yshiftu}=[shift = {(0, #1)}]

\tikzstyle{dot}=[inner sep=0.25mm,minimum width=1mm,minimum height=1mm,draw,shape=circle,text depth=-0.2mm]
\tikzstyle{white dot}=[dot,fill=white, draw=black]
\tikzstyle{action}=[dot,fill=white,scale=0.667,inner sep=0.5mm]
\tikzstyle{copier}=[dot,fill=white,scale=2.0]
\tikzstyle{black copier}=[dot,fill=black,scale=2.0]

\tikzstyle{box}=[fill=white, draw=black, shape=rectangle]
\tikzstyle{medium box}=[fill=white, draw=black, shape=rectangle, minimum width=1.5cm, minimum height=0.66cm]
\tikzstyle{arrow box}=[fill=white, draw, shape=rectangle,minimum height=5mm,yshift=-0.5mm,minimum width=5mm]

\tikzstyle{effect}=[regular polygon, regular polygon sides=3,draw]
\tikzstyle{state0}=[regular polygon, regular polygon sides=3,draw,shape border rotate=0]
\tikzstyle{state90}=[regular polygon, regular polygon sides=3,draw,shape border rotate=90]
\tikzstyle{state180}=[regular polygon, regular polygon sides=3,draw,shape border rotate=180]
\tikzstyle{state270}=[regular polygon, regular polygon sides=3,draw,shape border rotate=270]

\tikzstyle{scalar}=[diamond,draw,inner sep=1pt]

\tikzstyle{discarder}=[my ground,draw,inner sep=0pt,minimum width=4.2pt,minimum height=11.2pt,anchor=input,rotate=90]
\tikzstyle{discarder0}=[my ground,draw,inner sep=0pt,minimum width=4.2pt,minimum height=11.2pt,anchor=input,rotate=0]

\tikzstyle{pointy1}=[->]
\tikzstyle{midpoint1}=[-, {postaction={decorate,decoration={markings, mark=at position .5 with {\arrow{>}}}}}]
\tikzstyle{midpointy1pointy}=[->, {postaction={decorate,decoration={markings, mark=at position .5 with {\arrow{>}}}}}]
\tikzstyle{dashed1}=[-, dashed]
\tikzstyle{dotted1}=[-, dotted]
\tikzstyle{dash-pointy}=[->, dashed]

\input{strings.tikzdefs}

\ifexternalizetikz\tikzexternaldisable\fi
\newsavebox\sbground
\savebox\sbground{%
  \begin{tikzpicture}[baseline=0pt]
    \draw (0,-.1ex) to (0,.85ex)
    node[ground IEC,draw,anchor=input,inner sep=0pt,
    minimum width=3.15pt,minimum height=8.4pt,rotate=90] {};
  \end{tikzpicture}%
}
\newcommand{\ground}{\mathord{\usebox\sbground}}

\newsavebox\sbcopier
\savebox\sbcopier{%
  \begin{tikzpicture}[baseline=0pt]
    \node[copier,scale=0.7] (a) at (0,3.8pt) {};
    \draw (a) -- +(-90:.21);
    \draw (a) -- +(45:.21);
    \draw (a) -- +(135:.21);
  \end{tikzpicture}}

\ifexternalizetikz\tikzexternalenable\fi

\newsavebox\bsbcopier
\savebox\bsbcopier{%
  \begin{tikzpicture}[baseline=0pt]
    \node[black copier,scale=0.7] (a) at (0,3.8pt) {};
    \draw (a) -- +(-90:.21);
    \draw (a) -- +(45:.21);
    \draw (a) -- +(135:.21);
  \end{tikzpicture}}

\ifexternalizetikz\tikzexternalenable\fi

\usepackage{quiver}

\date{\today}

\newcommand{\Poly}[1]{\Cat{Poly}_{\cat{#1}}}
\newcommand{\Dyn}{\Cat{Dyn}}
\newcommand{\DynT}[1]{\Dyn^{\mathbb{#1}}}
\newcommand{\RDyn}{\Cat{RDyn}}
\newcommand{\RDynT}[1]{\RDyn^{\mathbb{#1}}}
\newcommand{\BDyn}{\Cat{BunDyn}}
\newcommand{\BDynT}[1]{\BDyn^{\mathbb{#1}}}
\newcommand{\RBDyn}{\Cat{RBDyn}}
\newcommand{\RBDynT}[1]{\RBDyn^{\mathbb{#1}}}
\newcommand{\NDyn}{\Cat{NDyn}}
\newcommand{\NDynT}[1]{\NDyn^{\mathbb{#1}}}
\newcommand{\RNDyn}{\Cat{RNDyn}}
\newcommand{\RNDynT}[1]{\RNDyn^{\mathbb{#1}}}
\newcommand{\pTCoalg}[2]{{#1}{#2}\mdash\Cat{Coalg}}
\newcommand{\pTCoalgT}[3]{\pTCoalg{#1}{#2}^{#3}}

\newcommand{\PTCoalgT}[1]{\pTCoalgT{p}{T}{#1}}
\newcommand{\TCoalgT}[1]{\pTCoalgT{(-)}{T}{#1}}

\newcommand{\deloop}[1]{\mathbf{B}#1}
\newcommand{\Sum}{\sum\limits}

\author{Toby St. Clere Smithe}
\date{\today}
\title{Some Notions of (Open) Dynamical System on Polynomial Interfaces}
\begin{document}

\maketitle
\begin{abstract}
  We define indexed categories of (open) dynamical system and random dynamical
  system over polynomial interfaces, where time is given by an arbitrary monoid
  \(\Tt\). We consider the case of open random dynamical systems over both open
  and closed noise sources, and the case where the interface of the random
  system is `nested' over the interface of its noise source. We show that, in
  discrete time, our categories of dynamical systems over polynomial interfaces
  \(p\) are equivalent to Spivak's categories \(p\mdash\Cat{Coalg}\) of
  \(p\)-coalgebras. We then define a notion of generalized \(pT\)-coalgebra for
  a monad \(T\), thereby extending the coalgebraic notion of dynamical system to
  general time, and show that this construction bestows a notion of open Markov
  process when the monad \(T\) is a probability monad. Finally, we speculate on
  some further connections and open questions.
\end{abstract}

\section{Polynomials for embodiment and interaction}
\label{sec:org2c00a4c}

Each system in our universe inhabits some interface or boundary. It receives
signals from its environments through this boundary, and can act by changing its
shape (and, as we will see later, its position). As a system changes its shape,
the set of possible immanent signals might change accordingly: consider a
hedgehog rolling itself into a ball, thereby protecting its soft underbelly from
harm (amongst other immanent signals). A system may also change its shape by
coupling itself to some other system, such as when we pick up chalk to work
through a problem. And shapes can be abstract: we change our `shapes' when we
enter an online video conference, or move within a virtual reality. We describe
all of these interactions formally using polynomial functors, drawing on the
work of \textcite{Spivak2020Poly}.

\begin{defn}
  Let $\cat{E}$ be a locally Cartesian closed category, and denote by $y^A$ the
  representable copresheaf $y^A := \cat{E}(A, -) : \cat{E} \to \cat{E}$. A
  \emph{polynomial functor} $p$ is a coproduct of representable functors,
  written $p := \sum_{i : p(1)} y^{p_i}$, where $p(1) : \cat{E}$ is the indexing
  object. The category of polynomial functors in $\cat{E}$ is the full
  subcategory $\Poly{E} \hookrightarrow [\cat{E}, \cat{E}]$ of the
  \(\cat{E}\)-copresheaf category spanned by coproducts of representables. A
  morphism of polynomials is therefore a natural transformation.
\end{defn}

\begin{rmk} %
  Every copresheaf $P : \cat{E} \to \cat{E}$ corresponds to a bundle $p : E \to
  B$ in $\cat{E}$, for which $B = P(1)$ and for each $i : P(1)$, the fibre $p_i$
  is $P(i)$. We will henceforth elide the distinction between a copresheaf $P$
  and its corresponding bundle $p$, writing $p(1) := B$ and $p[i] := p_i$, where
  $E = \sum_i p[i]$. A natural transformation $f : p \to q$ between copresheaves
  therefore corresponds to a map of bundles. In the case of polynomials, by the
  Yoneda lemma, this map is given by a `forwards' map $f_1 : p(1) \to q(1)$ and
  a family of `backwards' maps $f^\# : q[f_1(\mdash)] \to p[\mdash]$ indexed by
  $p(1)$, as in the left diagram below. Given $f : p \to q$ and $g : q \to r$,
  their composite $g \circ f : p \to r$ is as in the right diagram below.
  \begin{equation*}
    \begin{tikzcd}
      E & {f^*F} & F \\
      B & B & C
      \arrow["{f^\#}"', from=1-2, to=1-1]
      \arrow[from=1-2, to=1-3]
      \arrow["q", from=1-3, to=2-3]
      \arrow["p"', from=1-1, to=2-1]
      \arrow[from=2-1, to=2-2, Rightarrow, no head]
      \arrow["{f_1}", from=2-2, to=2-3]
      \arrow[from=1-2, to=2-2]
      \arrow["\lrcorner"{anchor=center, pos=0.125}, draw=none, from=1-2, to=2-3]
    \end{tikzcd}
    \qquad\qquad
    \begin{tikzcd}
      E & {f^*g^*G} & G \\
      B & B & D
      \arrow["{(gf)^\#}"', from=1-2, to=1-1]
      \arrow[from=1-2, to=1-3]
      \arrow["r", from=1-3, to=2-3]
      \arrow["p"', from=1-1, to=2-1]
      \arrow[from=2-2, to=2-1, Rightarrow, no head]
      \arrow["{g_1 \circ f_1}", from=2-2, to=2-3]
      \arrow[from=1-2, to=2-2]
      \arrow["\lrcorner"{anchor=center, pos=0.125}, draw=none, from=1-2, to=2-3]
    \end{tikzcd}
  \end{equation*}
  where $(gf)^\#$ is given by the \(p(1)\)-indexed family of composite maps
  $r[g_1(f_1(\mdash))] \xto{f^\ast g^\#} q[f_1(\mdash)] \xto{f^\#} p[\mdash]$.
\end{rmk}

In our morphological semantics, we will call a polynomial \(p\) a
\emph{phenotype}, its base type \(p(1)\) its \emph{morphology} and the total
space \(\sum_i p[i]\) its \emph{sensorium}. We will call elements of the
morphology \emph{shapes} or \emph{configurations}, and elements of the sensorium
\emph{immanent signals}.

\begin{prop}[\textcite{Spivak2020Poly}]
  There is a monoidal structure $(\Poly{E}, \otimes, y)$ that we interpret as
  ``putting systems in parallel''. Given $p : \sum_i p[i] \to p(1)$ and $q :
  \sum_j q[j] \to q(1)$, we have $p \otimes q = \sum_i \sum_j p[i] \times q[j]
  \to p(1) \times q(1)$. $y : 1 \to 1$ is then clearly unital. \qed
\end{prop}

\begin{prop}[\textcite{Spivak2020Poly}]
  The monoidal structure $(\Poly{E}, \otimes, y)$ is closed, with corresponding
  internal hom denoted $[-, -]$.
\end{prop}

We interpret morphisms \((f_1, f^\#)\) of polynomials as encoding interaction
patterns; in particular, such morphisms encode how composite systems act as
unities. For example, a morphism \(f : p \otimes q \to r\) specifies how the
systems \(p\) and \(q\) come together to form a system \(r\): the map \(f_1\)
encodes how \(r\)-configurations are constructed from configurations of \(p\)
and \(q\); and the map \(f^\#\) encodes how immanent signals on \(p\) and \(q\)
result from signals on \(r\) or from the interaction of \(p\) and \(q\). For
intuition, consider two people engaging in a handshake, or an enzyme acting on a
protein to form a complex. The internal hom \([o, p]\) encodes all the possible
ways that an \(o\)-phenotype system can ``plug into'' a \(p\)-phenotype system.

\begin{rmk}
  In the literature on active inference and the free energy principle, there is
  much debate about the concept of `Markov blanket', an informal notion
  conceived to represent the boundary of an adaptive system. We believe that the
  algebra of polynomials is sufficient to formalize this concept precisely, and
  clear up much of the confusion in the literature.
\end{rmk}

\section{Nested systems and dependent polynomials}
\label{sec:org3cddd26}

The polynomial formalism as presented in the previous section suffices to
describe systems' shapes, and behaviours of those shapes that depend on their
sensoria. But in our world, a system has a \emph{position} as well as a shape!
Indeed, one might want to consider systems nested within systems, such that the
outer systems constitute the `universes' of the inner systems; in this way,
inner shapes may depend on outer shapes, and inner sensoria on outer
sensoria.\footnote{We might even consider the outer shapes explicitly as
positions in some world-space, and the outer sensorium as determined by possible
paths between positions, in agreement with the perspective of
\textcite{Spivak2020Poly} on polynomials.} We can model this situation
polynomially.

Recall that an object in \(\Poly{E}\) corresponds to a bundle \(E \to B\),
equivalently a diagram \(1 \leftarrow E \xto{p} B \to 1\), and note that the
unit polynomial \(y\) corresponds to a bundle \(1 \to 1\). We can then think of
\(\Poly{E}\) as the category of ``polynomials in one variable'', or
``polynomials over y''. This presents a natural generalization, to polynomials
in many variables, corresponding to diagrams \(J \leftarrow E \to B \to I\);
these diagrams form the objects of a category \(\Poly{E}(J, I)\). When \(J\) is
a (polynomial) bundle \(\beta\) over \(I\), then we can take the subcategory of
\(\Poly{E}(J, I)\) whose objects are commuting squares and whose morphisms are
prisms as follows; the commutativity ensures that inner and outer sensoria are
compatible.

\begin{prop} \label{prop:poly-nest}
  There is an indexed category of \textbf{nested polynomials} which by abuse of
  notation we will call $\Poly{E}(-) : \Poly{E} \to \Cat{Cat}$. Given $\beta : J
  \to I$, the category $\Poly{E}(\beta)$ has commuting squares as on the left
  below as objects and prisms as on the right as morphisms. Its action on
  polynomial morphisms $\beta \to \gamma$ is given by composition.
  
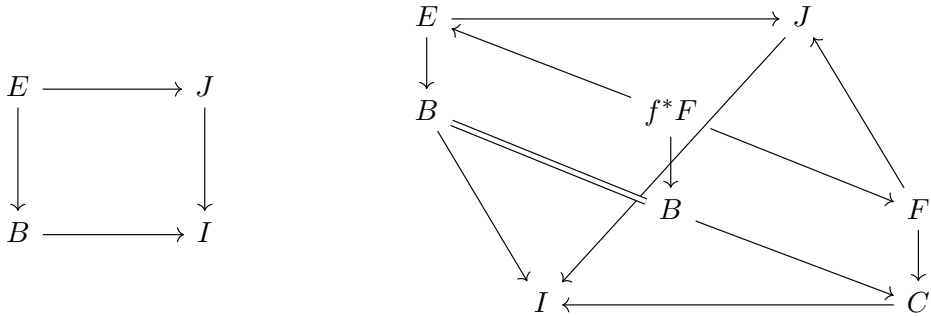
\begin{equation*}
    \begin{tikzcd}
      E && J \\
      \\
      B && I
      \arrow[from=1-3, to=3-3]
      \arrow[from=1-1, to=1-3]
      \arrow[from=1-1, to=3-1]
      \arrow[from=3-1, to=3-3]
    \end{tikzcd}
    \qquad\qquad\qquad
    \begin{tikzcd}
      E &&& J \\
      B && {f^*F} \\
      && B && F \\
      & I &&& C
      \arrow[from=1-4, to=4-2]
      \arrow[from=2-1, to=4-2]
      \arrow[Rightarrow, no head, from=2-1, to=3-3]
      \arrow[from=3-3, to=4-5]
      \arrow[from=2-3, to=3-5]
      \arrow[from=3-5, to=4-5]
      \arrow[from=2-3, to=1-1]
      \arrow[from=1-1, to=2-1]
      \arrow[from=4-5, to=4-2]
      \arrow[from=3-5, to=1-4]
      \arrow[from=1-1, to=1-4]
      \arrow[from=2-3, to=3-3]
    \end{tikzcd}
  \end{equation*}
\end{prop}

\begin{prop}
  The base category of polynomials \(\Poly{E}\) is isomorphic to its image over
  the trivial polynomial \(y\): \(\Poly{E} \cong \Poly{E}(y)\).
\end{prop}

\begin{rmk}
  This construction can be repeatedly iterated, modelling systems within systems
  within systems. We leave the consideration of the structure of this iteration
  to future work, though we expect it to have a ``mutually coinductive'' type
  and an opetopic shape, perhaps equivalent to that obtained by iterating the
  $\Cat{Para}$ construction.
\end{rmk}

\section{Dynamical systems with polynomial interfaces}
\label{sec:orge301230}

We describe categories of dynamical systems. We begin with classical `closed'
deterministic systems, before defining closed measure-preserving systems and
closed random dynamical systems: we can think of random dynamical systems as
deterministic systems parameterized by some dynamical noise source. More
precisely, a random dynamical system will be a bundle of a dynamical system over
a measure-preserving dynamical system on some probability space.

We then progressively open these various system types up to polynomial
interaction: first, general (deterministic) open dynamical systems; proceeded by
open random dynamical systems over closed bases, open measure-preserving
systems, and then ``fully open'' random dynamical systems (with open
measure-preserving bases). Finally, we consider the case where an open random
dynamical system is defined on a nested polynomial interface: that is, where the
interface of the total random system is nested in the interface of the base
measure-preserving system. This formalizes the idea that randomness comes from
some noise source that is ``really external'' to the random system, in such a
way that it pervades the random system's environment: think for instance of the
cosmic microwave background in our own universe. In order for this randomness to
be compatible with any change-of-scale implied by the nesting of open systems,
we also require that the dynamics be accordingly compatible. At the end of the
section, we consider the connection between our categories of dynamical systems
an \(p\)-coalgebras defined internally to \(\Poly{E}\).

Now, however, our starting point is classical closed systems.

\begin{defn}
  Let \((\Tt, +, 0)\) be a monoid, representing time. Let \(X : \cat{E}\) be
  some space, called the \textbf{state space}.  Then a \textbf{closed dynamical
    system} \(\vartheta\) \textbf{with state space} \(X\) \textbf{and time}
  \(\Tt\) is an action of \(\Tt\) on \(X\) satisfying a \textit{flow condition}.
  When \(\Tt\) is also an object of \(\cat{E}\), then this amounts to a morphism
  \(\vartheta : \Tt \times X \to X\) (or equivalently, a time-indexed family of
  \(X\)-endomorphisms, \(\vartheta(t) : X \to X\)), such that \(\vartheta(0) =
  \id_X\) and \(\vartheta(s + t) = \vartheta(s) \circ \vartheta(t)\). We will
  call these criteria the \textbf{flow condition} for \(\vartheta\).
\end{defn}

\begin{prop} \label{prop:transition-map}
  When time is discrete, as in the case \(\Tt = \nn\), any dynamical system
  \(\vartheta\) is entirely determined by its action at \(1 : \Tt\). That is,
  letting the state space be \(X\), we have \(\vartheta(t) = \vartheta(1)^{\circ
    t}\) where \(\vartheta(1)^{\circ t}\) means ``compose \(\vartheta(1) : X \to
  X\) with itself \(t\) times''.
  \begin{proof}
    The proof is by (co)induction on \(t : \Tt\). We must have \(\vartheta(0) =
    \id_X\) and \(\vartheta(t+s) = \vartheta(t) \circ \vartheta(s)\). So for any
    \(t\), we must have \(\vartheta(t+1) = \vartheta(t) \circ
    \vartheta(1)\). The result follows immediately; note for example that
    \(\vartheta(2) = \vartheta(1+1) = \vartheta(1) \circ \vartheta(1)\).
  \end{proof}
\end{prop}

\begin{ex} \label{ex:closed-vector-field}
  Suppose \(X : M \to TM\) is a vector field on \(M\), with a corresponding
  solution (integral curve) \(\chi_x : \rr \to M\) for all \(x : M\); that is,
  \(\chi'(t) = X(\chi_x(t))\) and \(\chi_x(0) = x\). Then letting the point
  \(x\) vary, we obtain a map \(\chi : \rr \times M \to M\). This \(\chi\) is a
  closed dynamical system with state space \(M\) and time \(\rr\).
\end{ex}

\begin{prop} \label{prop:closed-dyn-cat}
  Closed dynamical systems with state spaces in \(\cat{E}\) and time \(\Tt\) are
  the objects of the functor category \(\Cat{Cat}(\deloop{\Tt}, \cat{E})\),
  where \(\deloop{\Tt}\) denotes the delooping of the monoid \(\Tt\).  Morphisms
  of dynamical systems are therefore natural transformations.
  \begin{proof}
    The category \(\deloop{\Tt}\) has a single object \(\ast\) and morphisms \(t
    : \ast \to \ast\) for each point \(t : \Tt\); the identity is the monoidal
    unit \(0 : \Tt\) and composition is given by \(+\).  A functor \(\vartheta :
    \deloop{\Tt} \to \cat{E}\) therefore picks out an object \(\vartheta(\ast) :
    \cat{E}\), and, for each \(t : \Tt\), a morphism \(\vartheta(t) :
    \vartheta(\ast) \to \vartheta(\ast)\), such that the functoriality condition
    is satisfied. Functoriality requires that identities map to identities and
    composition is preserved, so we require that \(\vartheta(0) =
    \id_{\vartheta(\ast)}\) and that \(\vartheta(s + t) = \vartheta(s) \circ
    \vartheta(t)\). Hence the data for a functor \(\vartheta : \deloop{\Tt} \to
    \cat{E}\) amount to the data for a closed dynamical system in \(\cat{E}\)
    with time \(\Tt\), and the functoriality condition amounts precisely to the
    flow condition. A morphism of closed dynamical systems \(f : \vartheta \to
    \psi\) is a map on the state spaces \(f : \vartheta(\ast) \to \psi(\ast)\)
    that commutes with the flow, meaning that \(f\) satisfies \(f \circ
    \vartheta(t) = \psi(t) \circ f\) for all times \(t : \Tt\); this is
    precisely the definition of a natural transformation \(f : \vartheta \to
    \psi\) between the corresponding functors.
  \end{proof}
\end{prop}

We now consider closed \emph{random} dynamical systems.

\begin{defn}
  Suppose \(\cat{E}\) is a category equipped with a probability monad \(\Pa :
  \cat{E} \to \cat{E}\) and a terminal object \(1 : \cat{E}\).  A
  \textbf{probability space} in \(\cat{E}\) is an object of the slice \(1 /
  \Kl(\Pa)\) of the Kleisli category of the probability monad under
  \(1\). Explicitly, a probability space is therefore equivalently a pair \((B,
  \beta)\) where \(B : \cat{E}\) is an object and \(\beta : 1 \to \Pa B\) is a
  measure over \(B\). Morphisms \(f : (A, \alpha) \to (B, \beta)\) between
  probability spaces are stochastic channels \(f : A \to \Pa B\) that preserve
  the measure; that is, they satisfy \(f \klcirc \alpha = \beta\).
\end{defn}

\begin{prop}
  There is a forgetful functor \(1/\Kl(\Pa) \to \cat{E}\) taking probability
  spaces \((B, \beta)\) to the underlying spaces \(B\) and their morphisms \(f :
  (A, \alpha) \to (B, \beta)\) to the underlying maps \(f : A \to \Pa B\). We
  will write \(B\) to refer to the space in \(\cat{E}\) underlying a probability
  space \((B, \beta)\), in the image of this forgetful functor.
\end{prop}

\begin{defn} \label{def:metric-sys}
  Let \((B, \beta)\) be a probability space in \(\cat{E}\). A closed
  \textbf{metric} or \textbf{measure-preserving} dynamical system \((\vartheta,
  \beta)\) on \((B, \beta)\) with time \(\Tt\) is a closed dynamical system
  \(\vartheta\) with state space \(B : \cat{E}\) such that, for all \(t : \Tt\),
  \(\Pa \vartheta(t) \circ \beta = \beta\); that is, each \(\vartheta(t)\) is a
  \((B, \beta)\)-endomorphism in \(1/\Kl(\Pa)\).
\end{defn}

\begin{prop} \label{prop:metric-sys}
  Closed measure-preserving dynamical systems in \(\cat{E}\) with time \(\Tt\)
  form the objects of a category \(\Cat{Cat}(\deloop{\Tt}, \cat{E})_{\Pa}\)
  whose morphisms \(f : (\vartheta, \alpha) \to (\psi, \beta)\) are maps \(f :
  \vartheta(\ast) \to \psi(\ast)\) in \(\cat{E}\) between the state spaces that
  preserve both flow and measure, as in the following commutative diagram, which
  also indicates their composition:
  \[\begin{tikzcd}
	&& {\Pa\vartheta(\ast)} && {\Pa\vartheta(\ast)} \\
	\\
	1 && {\Pa\psi(\ast)} && {\Pa\psi(\ast)} && 1 \\
	\\
	&& {\Pa\lambda(\ast)} && {\Pa\lambda(\ast)}
	\arrow["\alpha", from=3-1, to=1-3]
	\arrow["\beta"', from=3-1, to=3-3]
	\arrow["\gamma"', from=3-1, to=5-3]
	\arrow["\alpha"', from=3-7, to=1-5]
	\arrow["\beta", from=3-7, to=3-5]
	\arrow["\gamma", from=3-7, to=5-5]
	\arrow["{\Pa\vartheta(t)}", from=1-3, to=1-5]
	\arrow["{\Pa\psi(t)}"', from=3-3, to=3-5]
	\arrow["{\Pa\lambda(t)}"', from=5-3, to=5-5]
	\arrow["{\Pa f}", from=1-3, to=3-3]
	\arrow["{\Pa f}"', from=1-5, to=3-5]
	\arrow["{\Pa g}", from=3-3, to=5-3]
	\arrow["{\Pa g}"', from=3-5, to=5-5]
  \end{tikzcd}\]
  \begin{proof}
    The identity morphism on a closed measure-preserving dynamical system is the
    identity map on its state space.  It is easy to check that composition as in
    the diagram above is thus both associative and unital with respect to these
    identities.
  \end{proof}
\end{prop}

\begin{prop} \label{prop:forget-metric-sys}
  There is a forgetful functor \(U : \Cat{Cat}(\deloop{\Tt}, \cat{E})_{\Pa} \to
  \Cat{Cat}(\deloop{\Tt}, \cat{E})\) which simply forgets the probability space
  structures. Given a closed measure-preserving dynamical system \(\vartheta,
  \beta\), we will write \(\vartheta\) for the closed dynamical system in the
  image of this forgetful functor.
\end{prop}

\begin{defn}
  Let \((\vartheta, \beta)\) be a closed measure-preserving dynamical system. A
  closed \textbf{random dynamical system} over \((\vartheta, \beta)\) is an
  object of the slice category \(\Cat{Cat}(\deloop{\Tt}, \cat{E})/\vartheta\);
  it is therefore a bundle of the corresponding functors.
\end{defn}

\begin{prop}
  The indexing of categories of closed random dynamical systems
  \(\Cat{Cat}(\deloop{\Tt}, \cat{E})/\vartheta\) by their base
  measure-preserving systems \((\vartheta, \beta)\) constitutes an indexed
  category \(\Cat{Cat}(\deloop{\Tt}, \cat{E})_{\Pa}\op \to \Cat{Cat}\) taking
  closed measure-preserving systems to the categories of closed random dynamical
  systems (bundles of functors) above them, and morphisms of closed
  measure-preserving systems to the corresponding base-change functors. This
  indexed category is therefore defined as \(\Cat{Cat}(\deloop{\Tt},
  \cat{E})/U\op(-)\) where \(U\op : \Cat{Cat}(\deloop{\Tt}, \cat{E})_{\Pa}\op
  \to \Cat{Cat}(\deloop{\Tt}, \cat{E})\op\) is the opposite of the forgetful
  functor in Proposition \ref{prop:forget-metric-sys}.
\end{prop}

\begin{ex} \label{ex:brown-sde}
  The solutions \(X(t, \omega; x_0) : \rr_+ \times \Omega \times M \to M\) to a
  stochastic differential equation \(\d X_t = f(t, X_t) \d t + \sigma(t, X_t) \d
  W_t\), where \(W : \rr_+ \times \Omega \to M\) is a Wiener process in \(M\),
  define a random dynamical system \(\rr_+ \times \Omega \times M \to M : (t,
  \omega, x) \mapsto X(t, \omega; x_0)\) over the Wiener base flow \(\theta :
  \rr_+ \times \Omega \to \Omega : (t, \omega) \mapsto W(s+t, \omega) - W(t,
  \omega)\) for any \(s : \rr_+\). We can alternatively represent this system as
  a bundle system over \((\theta, \gamma)\), where \(\gamma\) is the Wiener
  measure on the Wiener space \(\Omega\), by writing \(\vartheta : \rr_+ \times
  \Omega \times M \to \Omega \times M : (t, \omega, x) \mapsto \left(\theta(t,
  \omega), X(t, \omega; x_0)\right)\). This gives a closed random dynamical
  system in \(\Cat{Cat}(\deloop{\rr_+}, \cat{E})/\theta\).
\end{ex}

\begin{ex} \label{ex:markov-chain}
  In discrete time (\(\Tt = \nn\)), closed random dynamical systems (with
  independent-increment noise) correspond to Markov chains, for the same reason
  that general closed discrete-time dynamical systems correspond to transition
  functions (Proposition \ref{prop:transition-map}). Let \(\vartheta\) be a
  closed discrete-time random dynamical system over \((\theta, \gamma)\), and
  suppose the state space of \(\vartheta\) corresponds to a trivial bundle \(\pi
  : \Omega \times M \to \Omega\). Since the systems \(\vartheta\) and \(\theta\)
  correspond to transition functions, we have \(\vartheta^u : \Omega \times M
  \to \Omega \times M\) and \(\theta^u : \Omega \to \Omega\). By the universal
  property of the product, the map \(\vartheta^u : \Omega \times M \to \Omega
  \times M\) corresponds to a pair of maps \(\vartheta^\Omega : \Omega \times M
  \to \Omega\) and \(\vartheta^M : \Omega \times M \to M\), but the former
  component must coincide with \(\theta^u\), since \(\vartheta\) is a bundle
  over \(\theta\). At each time step \(n : \nn\), the noise is distributed
  according to \(\gamma : 1 \to \Pa \Omega\). Pushing \(\gamma\) forward along
  \(\vartheta^M\) induces a \(\Pa\)-coalgebra, \(\vartheta^\gamma : M \to \Pa
  M\), which is precisely a Markov chain. Conversely, note that, by randomness
  pushback \parencite[Def. 11.19]{Fritz2019synthetic}, any such map \(\tau : M
  \to \Pa M\) canonically induces a pair \((\tau^\flat : \Omega \times M \to M,
  \gamma : 1 \to \Pa \Omega)\) for some probability space \((\Omega, \gamma)\),
  and one can construct a random dynamical system accordingly; see
  \textcite[Theorem 2.1.6]{Arnold1998Random}.
\end{ex}

\begin{defn} \label{def:poly-dyn}
  Let \(p : \Poly{E}\) be a polynomial in \(\cat{E}\), to be called the
  \textbf{interface}.  Let \(S : \cat{E}\) be an object, to be called the
  \textbf{state space}.  Let \((\Tt : \cat{E}, +, 0)\) be a monoid, representing
  time. An \textbf{open dynamical system on the interface} \(p\) \textbf{with
    state space} \(S\) \textbf{and time} \(\Tt\) consists in a pair of morphisms
  \(\vartheta^o : \Tt \times S \to p(1)\) and \(\vartheta^u : \Sum_{t : \Tt}
  \Sum_{s : S} p[\vartheta^o(t, s)] \to S\), such that, for any global section
  \(\sigma : p(1) \to \Sum_{i:p(1)} p[i]\) of \(p\), the maps \(\vartheta^\sigma
  : \Tt \times S \to S\) given by
  \[
  \Sum_{t:\Tt} S \xto{\vartheta^o(-)^\ast \sigma} \Sum_{t:\Tt} \Sum_{s:S} p[\vartheta^o(-, s)] \xto{\vartheta^u} S
  \]
  constitute a \textit{closed} dynamical system, \textit{i.e.}, an object in
  \(\Cat{Cat}(\deloop{\Tt}, \cat{E})\). That is, the maps \(\vartheta^\sigma\)
  must satisfy the flow conditions, that \(\vartheta^\sigma(0) = \id_S\) and
  \(\vartheta^\sigma(s + t) = \vartheta^\sigma(s) \circ
  \vartheta^\sigma(t)\). We collect the data of such a dynamical system over
  \(p\) into a tuple \(\vartheta = (S, \vartheta^o, \vartheta^u)\). We call the
  closed system \(\vartheta^\sigma\), induced by a section \(\sigma\) of \(p\),
  the \textbf{closure} of \(\vartheta\) by \(\sigma\).
\end{defn}

\begin{prop} \label{prop:poly-dyn}
  Open dynamical systems over \(p\) with time \(\Tt\) form a category, denoted
  \(\DynT{T}(p)\).  Its morphisms are defined as follows.  Let \(\vartheta :=
  (X, \vartheta^o, \vartheta^u)\) and \(\psi := (Y, \psi^o, \psi^u)\) be two
  dynamical systems over \(p\). A morphism \(f : \vartheta \to \psi\) consists
  in a morphism \(f : X \to Y\) such that, for any time \(t : \Tt\) and global
  section \(\sigma : p(1) \to \Sum_{i:p(1)} p[i]\) of \(p\), the following
  naturality squares commute:
  \[\begin{tikzcd}
	X & {\Sum_{x:X} p[\vartheta^o(t, x)]} & X \\
	Y & {\Sum_{y:Y} p[\psi^o(t, y)]} & Y
	\arrow["{\vartheta^o(t)^\ast \sigma}", from=1-1, to=1-2]
	\arrow["{\vartheta^u(t)}", from=1-2, to=1-3]
	\arrow["f"', from=1-1, to=2-1]
	\arrow["f", from=1-3, to=2-3]
	\arrow["{\psi^o(t)^\ast \sigma}"', from=2-1, to=2-2]
	\arrow["{\psi^u(t)}"', from=2-2, to=2-3]
  \end{tikzcd}\]
  The identity morphism \(\id_\vartheta\) on the dynamical system \(\vartheta\)
  is given by the identity morphism \(\id_X\) on its state space
  \(X\). Composition of morphisms of dynamical systems is given by composition
  of the morphisms of the state spaces.
  \begin{proof}
    We need to check unitality and associativity of composition. This amounts to
    checking that the composite naturality squares commute. But this follows
    immediately, since the composite of two commutative diagrams along a common
    edge is again a commutative diagram.
  \end{proof}
\end{prop}

\begin{prop} \label{prop:poly-dyn-idx}
  \(\DynT{T}(p)\) extends to a polynomially-indexed category, \(\DynT{T} :
  \Poly{E} \to \Cat{Cat}\). Suppose \(\varphi : p \to q\) is a morphism of
  polynomials.  We define a corresponding functor \(\DynT{T}(\varphi) :
  \DynT{T}(p) \to \DynT{T}(q)\) as follows.  Suppose \((X, \vartheta^o,
  \vartheta^u) : \DynT{T}(p)\) is an object (dynamical system) in
  \(\DynT{T}(p)\). Then \(\DynT{T}(\varphi)(X, \vartheta^o, \vartheta^u)\) is
  defined as the triple \((X, \varphi_1 \circ \vartheta^o, \vartheta^u \circ
  {\vartheta^o}^\ast \varphi^\#) : \DynT{T}(q)\), where the two maps are
  explicitly the following composites:
  \begin{gather*}
    \Tt \times X \xto{\vartheta^o} p(1) \xto{\varphi_1} q(1) \, ,
    \qquad
    \Sum_{t:\Tt} \Sum_{x:X} q[\varphi_1 \circ \vartheta^o(t, x)] \xto{{\vartheta^o}^\ast \varphi^\#} \Sum_{t:\Tt} \Sum_{x:X} p[\vartheta^o(t, x)] \xto{\vartheta^u} X \, .
  \end{gather*}
  On morphisms, \(\DynT{T}(\varphi)(f) : \DynT{T}(\varphi)(X, \vartheta^o,
  \vartheta^u) \to \DynT{T}(\varphi)(Y, \psi^o, \psi^u)\) is given by the same
  underlying map \(f : X \to Y\) of state spaces.
  \begin{proof}
    We need to check that \(\DynT{T}(\varphi)(X, \vartheta^o, \vartheta^u)\)
    satisfies the flow conditions of Definition \ref{def:poly-dyn}, that
    \(\DynT{T}(\varphi)(f)\) satisfies the naturality condition of Proposition
    \ref{prop:poly-dyn}, and that \(\DynT{T}\) is functorial with respect to
    polynomials. We begin with the flow condition. Given a section \(\tau : q(1)
    \to \Sum_{j:q(1)} q[j]\) of \(q\), we require the closures
    \(\DynT{T}(\varphi)(\vartheta)^\tau : \Tt \times X \to X\) given by
    \[
    \Sum_{t:\Tt} X \xto{\vartheta^o(-)^\ast \tau} \Sum_{t:\Tt} \Sum_{x:X} q[\varphi_1 \circ \vartheta^o(t, x)] \xto{{\vartheta^o}^\ast \varphi^\#} \Sum_{t:\Tt} \Sum_{x:X} p[\vartheta^o(t, x)] \xto{\vartheta^u} X
    \]
    to satisfy \(\DynT{T}(\varphi)(\vartheta)^\tau (0) = \id_X\) and
    \(\DynT{T}(\varphi)(\vartheta)^\tau (s+t) =
    \DynT{T}(\varphi)(\vartheta)^\tau (s) \circ
    \DynT{T}(\varphi)(\vartheta)^\tau (t)\). Note that the following diagram
    commutes, by the definition of \(\varphi^\#\),
    \[\begin{tikzcd}
	{\Sum_{i:p(1)} p[i]} && {\Sum_{i:p(1)}q[\varphi_1(i)]} && {p(1)} \\
	\\
	{p(1)} && {p(1)}
	\arrow["{\varphi_1^\ast q}"', from=1-3, to=3-3]
	\arrow[Rightarrow, no head, from=3-3, to=1-5]
	\arrow["{\varphi_1^\ast \tau}"', from=1-5, to=1-3]
	\arrow["{\varphi^\#}"', from=1-3, to=1-1]
	\arrow[Rightarrow, no head, from=3-1, to=3-3]
	\arrow["p"', from=1-1, to=3-1]
    \end{tikzcd}\]
    so that \(\varphi^\# \circ \varphi_1^\ast \tau\) is a section of
    \(p\). Therefore, letting \(\sigma := \varphi^\# \circ \varphi_1^\ast
    \tau\), for \(\DynT{T}(\varphi)(\vartheta)^\tau\) to satisfy the flow
    condition for \(\tau\) reduces to \(\vartheta^\sigma\) satisfying the flow
    condition for \(\sigma\). But this is given \textit{ex hypothesi} by
    Definition \ref{def:poly-dyn}, for any such section \(\sigma\), so
    \(\DynT{T}(\varphi)(\vartheta)^\tau\) satisfies the flow condition for
    \(\tau\). And since \(\tau\) was any section, we see that
    \(\DynT{T}(\varphi)(\vartheta)\) satisfies the flow condition generally.

    The proof that \(\DynT{T}(\varphi)(f)\) satisfies the naturality condition
    of Proposition \ref{prop:poly-dyn} proceeds similarly.  Supposing again that
    \(\tau\) is any section of \(q\), we require the following diagram to
    commute for any time \(t : \Tt\):
    \[\begin{tikzcd}
	X && {\Sum_{x:X} q[\varphi_1 \circ \vartheta^o(t, x)]} && {\Sum_{x:X} p[\vartheta^o(t, x)]} && X \\
	\\
	Y && {\Sum_{y:Y} q[\varphi_1 \circ \psi^o(t, x)]} && {\Sum_{y:Y} p[\psi^o(t, x)]} && Y
	\arrow["f"', from=1-1, to=3-1]
	\arrow["f", from=1-7, to=3-7]
	\arrow["{\vartheta^o(t)^\ast \varphi_1^\ast \tau}", from=1-1, to=1-3]
	\arrow["{\vartheta^o(t)^\ast \varphi^\#}", from=1-3, to=1-5]
	\arrow["{\vartheta^u(t)}", from=1-5, to=1-7]
	\arrow["{\psi^o(t)^\ast \varphi_1^\ast \tau}", from=3-1, to=3-3]
	\arrow["{\vartheta^o(t)^\ast \varphi^\#}", from=3-3, to=3-5]
	\arrow["{\psi^u(t)}", from=3-5, to=3-7]
    \end{tikzcd}\]
    Again letting \(\sigma := \varphi^\# \circ \varphi_1^\ast \tau\), we see
    that this diagram reduces exactly to the diagram in Proposition
    \ref{prop:poly-dyn} by the functoriality of pullback, and since \(f\) makes
    that diagram commute, it must also make this diagram commute.

    Finally, to show that \(\DynT{T}\) is functorial with respect to polynomials
    amounts to checking that composition and pullback are functorial; but this
    is a basic result of category theory.
  \end{proof}
\end{prop}

To confirm that our definition of open dynamical system subsumes the classical
case of `closed' dynamical systems, we now consider the case of dynamical
systems on the trivial interface \(y\): such a system has a trivial shape,
exposing no configuration to any environment nor receiving any signals from it.

\begin{prop} \label{prop:closed-sys-in-dyn-y}
  \(\DynT{T}(y)\) is equivalent to the classical category
  \(\Cat{Cat}(\deloop{\Tt}, \cat{E})\) of closed dynamical systems in
  \(\cat{E}\) with time \(\Tt\) (\textit{cf.} Proposition
  \ref{prop:closed-dyn-cat}).
  \begin{proof}
    The trivial interface \(y\) corresponds to the trivial bundle \(\id_1 : 1
    \to 1\). Therefore, a dynamical system over \(y\) consists of a choice of
    state space \(S\) along with a trivial output map \(\vartheta^o = \ground :
    \Tt \times S \to 1\) and a time-indexed update map \(\vartheta^u : \Tt
    \times S \to S\). This therefore has the form of a classical closed
    dynamical system, so it remains to check the flow condition. There is only
    one section of \(\id_1\), which is again \(\id_1\). Pulling this back along
    the unique map \(\vartheta^o(t) : S \to 1\) gives \(\vartheta^o(t)^\ast
    \id_1 = \id_S\). Therefore the requirement that, given any section
    \(\sigma\) of \(y\), each map \(\vartheta^u \circ \vartheta^o(t)^\ast
    \sigma\) satisfies the flow condition reduces to the classical requirement
    that \(\vartheta^u : \Tt \times S \to S\) satisfies the flow
    condition. Since the pullback of the unique section \(\id_1\) along the
    trivial output map \(\vartheta^o(t) = \ground : S \to 1\) of any dynamical
    system in \(\DynT{T}(y)\) is the identity of the corresponding state space
    \(\id_S\), a morphism \(f : (\vartheta(\ast), \vartheta^u, \ground) \to
    (\psi(\ast), \psi^u, \ground)\) in \(\DynT{T}(y)\) amounts precisely to a
    map \(f : \vartheta(\ast) \to \psi(\ast)\) on the state spaces in
    \(\cat{E}\) such that the naturality condition \(f \circ \vartheta^u(t) =
    \psi^u(t) \circ f\) of Proposition \ref{prop:closed-dyn-cat} is satisfied,
    and every morphism in \(\Cat{Cat}(\deloop{\Tt}, \cat{E})\) corresponds to a
    morphism in \(\DynT{T}(y)\) in this way.
  \end{proof}
\end{prop}

\begin{ex}
  We now consider the case of a dynamical system \((S, \vartheta^o,
  \vartheta^u)\) with outputs but no inputs.  Such a system lives over a
  polynomial bundle \(p : p(1) \xto{\sim} p(1)\) that is an isomorphism.  A
  section of this bundle must therefore be its inverse \(p^{-1}\), and so,
  \(\vartheta^o(t)^\ast p^{-1} = \id_{p(1)}\). Once again, the update map
  corresponds to a dynamical system in \(\Cat{Cat}(\deloop{\Tt}, \cat{E})\);
  just now we have outputs \(\vartheta^o : \Tt \times S \to p(1)\) exposed to
  the environment.
\end{ex}

\begin{prop} \label{prop:open-transition-map}
  When time is discrete, as with \(\Tt = \nn\), any open dynamical system \((X,
  \vartheta^o, \vartheta^u)\) over \(p\) is entirely determined by its
  components at \(1 : \Tt\). That is, we have \(\vartheta^o(t) = \vartheta^o(1)
  : X \to p(1)\) and \(\vartheta^u(t) = \vartheta^u(1) : \sum_{x:X}
  p[\vartheta^o(x)] \to X\). A discrete-time open dynamical system is therefore
  a triple \((X, \vartheta^o, \vartheta^u)\), where the two maps have types
  \(\vartheta^o : X \to p(1)\) and \(\vartheta^u : \sum_{x:X} p[\vartheta^o(x)]
  \to X\).
  \begin{proof}
    Suppose \(\sigma\) is a section of \(p\). We require each closure
    \(\vartheta^\sigma\) to satisfy the flow conditions, that
    \(\vartheta^\sigma(0) = \id_X\) and \(\vartheta^\sigma(t+s) =
    \vartheta^\sigma(t) \circ \vartheta^\sigma(s)\). In particular, we must have
    \(\vartheta^\sigma(t+1) = \vartheta^\sigma(t) \circ
    \vartheta^\sigma(1)\). By induction, this means that we must have
    \(\vartheta^\sigma(t) = \vartheta^\sigma(1)^{\circ t}\) (compare Proposition
    \ref{prop:transition-map}). Therefore we must in general have
    \(\vartheta^o(t) = \vartheta^o(1)\) and \(\vartheta^u(t) = \vartheta^u(1)\).
  \end{proof}
\end{prop}

\begin{ex}
  Suppose \(\dot{x} = f(x, a)\) and \(b = g(x)\), with \(f : X \times A \to TX\)
  and \(g : X \to B\). Then, as for the `closed' vector fields of Example
  \ref{ex:closed-vector-field}, this induces an open dynamical system \((X, \int
  f, g) : \DynT{R}(By^A)\), where \(\int f : \rr \times X \times A \to X\)
  returns the \((X,A)\)-indexed solutions of \(f\).
\end{ex}

\begin{ex}
  The preceding example is easily extended to the case of a general polynomial
  interface. Suppose similarly that \(\dot{x} = f(x, a_x)\) and \(b = g(x)\),
  now with \(f : \sum_{x:X} p[g(x)] \to TX\) and \(g : X \to p(1)\). Then we
  obtain an open dynamical system \((X, \int f, g) : \DynT{R}(p)\), where now
  \(\int f : \rr \times \sum_{x:X} p[g(x)] \to X\) is the `update' and \(g : X
  \to p(1)\) the `output' map.
\end{ex}

We now move on to define open \textit{random} dynamical systems. We do so in
stages, starting with defining open random dynamical systems over closed base
measure-preserving systems, and then opening up the base systems to form ``fully
open'' random dynamical systems. We then ask for the base and total open systems
to be compatible with appropriately nested corresponding polynomial interfaces,
resulting in a categorical notion of random nested dynamical system.

\begin{defn} \label{def:poly-rdyn}
  Let \((\theta, \beta)\) be a closed measure-preserving dynamical system in
  \(\cat{E}\) with time \(\Tt\), and let \(p : \Poly{E}\) be a polynomial in
  \(\cat{E}\). Write \(\Omega := \theta(\ast)\) for the state space of
  \(\theta\), and let \(\pi : S \to \Omega\) be an object (bundle) in
  \(\cat{E}/\Omega\). An \textbf{open random dynamical system over} \((\theta,
  \beta)\) \textbf{on the interface} \(p\) \textbf{with state space} \(\pi:S \to
  \Omega\) \textbf{and time} \(\Tt\) consists in a pair of morphisms
  \(\vartheta^o : \Tt \times S \to p(1)\) and \(\vartheta^u : \Sum_{t:\Tt}
  \Sum_{s:S} p[\vartheta^o(t, s)] \to S\), such that, for any global section
  \(\sigma : p(1) \to \Sum_{i:p(1)} p[i]\) of \(p\), the maps \(\vartheta^\sigma
  : \Tt \times S \to S\) defined as
  \[
  \Sum_{t:\Tt} S \xto{\vartheta^o(-)^\ast \sigma} \Sum_{t:\Tt} \Sum_{s:S} p[\vartheta^o(-, s)] \xto{\vartheta^u} S
  \]
  form a closed random dynamical system in \(\Cat{Cat}(\deloop{\Tt},
  \cat{E})/\theta\). That is to say, for all \(t : \Tt\) and sections
  \(\sigma\), the following naturality square commutes:
  \[\begin{tikzcd}
	S && {\Sum_{s:S} p[\vartheta^o(t, s)]} && S \\
	\Omega &&&& \Omega
	\arrow["\pi"', from=1-1, to=2-1]
	\arrow["\pi", from=1-5, to=2-5]
	\arrow["{\theta(t)}"', from=2-1, to=2-5]
	\arrow["{\vartheta^o(t)^\ast \sigma}", from=1-1, to=1-3]
	\arrow["{\vartheta^u(t)}", from=1-3, to=1-5]
  \end{tikzcd}\]
  We collect the data of such a dynamical system into a tuple \(\vartheta =
  (\pi, \vartheta^o, \vartheta^u)\). Given a section \(\sigma\) of \(p\), the
  induced closed system \(\vartheta^\sigma\) will again be called the
  \textbf{closure} of \(\vartheta\) by \(\sigma\).
\end{defn}

\begin{prop} \label{prop:poly-rdyn}
  Let \((\theta, \beta)\) be a closed measure-preserving dynamical system in
  \(\cat{E}\) with time \(\Tt\), and let \(p : \Poly{E}\) be a polynomial in
  \(\cat{E}\). Open random dynamical systems over \((\theta, \beta)\) on the
  interface \(p\) form the objects of a category \(\RDynT{T}(p,
  \theta)\). Writing \(\vartheta := (\pi_X, \vartheta^o, \vartheta^u)\) and
  \(\psi := (\pi_Y, \psi^o, \psi^u)\), a morphism \(f : \vartheta \to \psi\) is
  a map \(f: X \to Y\) in \(\cat{E}\) making the following diagram commute for
  all times \(t : \Tt\) and sections \(\sigma\) of \(p\):
  \[\begin{tikzcd}
	X &&&& {\Sum_{x:X} p[\vartheta^o(t, x)]} &&&& X \\
	\\
	&& \Omega &&&& \Omega \\
	\\
	Y &&&& {\Sum_{y:Y} p[\psi^o(t, y)]} &&&& Y
	\arrow["{\pi_X}"', from=1-1, to=3-3]
	\arrow["{\pi_X}", from=1-9, to=3-7]
	\arrow["{\theta(t)}", from=3-3, to=3-7]
	\arrow["{\vartheta^o(t)^\ast \sigma}", from=1-1, to=1-5]
	\arrow["{\vartheta^u(t)}", from=1-5, to=1-9]
	\arrow["{\psi^o(t)^\ast \sigma}"', from=5-1, to=5-5]
	\arrow["{\psi^u(t)}"', from=5-5, to=5-9]
	\arrow["{\pi_Y}", from=5-1, to=3-3]
	\arrow["{\pi_Y}"', from=5-9, to=3-7]
	\arrow["f"', from=1-1, to=5-1]
	\arrow["f", from=1-9, to=5-9]
  \end{tikzcd}\]
  Identities are given by the identity maps on state-spaces. Composition is
  given by pasting of diagrams.
\end{prop}

\begin{prop} \label{prop:poly-rdyn-idx}
  The categories \(\RDynT{T}(p, \theta)\) collect into a doubly-indexed category
  of the form \(\RDynT{T} : \Poly{E} \times \Cat{Cat}(\deloop{\Tt},
  \cat{E})_{\Pa} \to \Cat{Cat}\). By the universal property of the product
  \(\times\) in \(\Cat{Cat}\), it suffices to define the actions of
  \(\RDynT{T}\) separately on morphisms of polynomials and on morphisms of
  closed measure-preserving systems.

  Suppose therefore that \(\varphi : p \to q\) is a morphism of
  polynomials. Then, for each measure-preserving system \((\theta, \beta) :
  \Cat{Cat}(\deloop{\Tt}, \cat{E})_{\Pa}\), we define the functor
  \(\RDynT{T}(\varphi, \theta) : \RDynT{T}(p, \theta) \to \RDynT{T}(q, \theta)\)
  as follows. Let \(\vartheta := (\pi_X : X \to \Omega, \vartheta^o,
  \vartheta^u) : \RDynT{T}(p, \theta)\) be an object (open random dynamical
  system) in \(\RDynT{T}(p, \theta)\). Then, as in Proposition
  \ref{prop:poly-dyn-idx}, \(\RDynT{T}(\varphi, \theta)(\vartheta)\) is defined
  as the triple \((\pi_X, \varphi_1 \circ \vartheta^o, \vartheta^u \circ
  {\varphi^o}^\ast \varphi^\#) : \RDynT{T}(q, \theta)\), where the two maps are
  explicitly the following composites:
  \begin{gather*}
    \Tt \times X \xto{\vartheta^o} p(1) \xto{\varphi_1} q(1) \, ,
    \qquad
    \Sum_{t:\Tt} \Sum_{x:X} q[\varphi_1 \circ \vartheta^o(t, x)] \xto{{\vartheta^o}^\ast \varphi^\#} \Sum_{t:\Tt} \Sum_{x:X} p[\vartheta^o(t, x)] \xto{\vartheta^u} X \, .
  \end{gather*}
  Again as in Proposition \ref{prop:poly-dyn-idx}, on morphisms \(f : (\pi_X : X
  \to \Omega, \vartheta^o, \vartheta^u) \to (\pi_Y : Y \to \Omega, \psi^o,
  \psi^u)\), the image \(\RDynT{T}(\varphi, \theta)(f) : \RDynT{T}(\varphi,
  \theta)(\pi_X, \vartheta^o, \vartheta^u) \to \RDynT{T}(\varphi, \theta)(\pi_Y,
  \psi^o, \psi^u)\) is given by the same underlying map \(f : X \to Y\) of state
  spaces.

  Next, suppose that \(\phi : (\theta, \beta) \to (\theta', \beta')\) is a
  morphism of closed measure-preserving dynamical systems, and let \(\Omega' :=
  \theta'(\ast)\) be the state space of the system \(\theta'\). By Proposition
  \ref{prop:metric-sys}, the morphism \(\phi\) corresponds to a map \(\phi :
  \Omega \to \Omega'\) on the state spaces that preserves both flow and
  measure. Therefore, for each polynomial \(p : \Poly{E}\), we define the
  functor \(\RDynT{T}(p, \phi) : \RDynT{T}(p, \theta) \to \RDynT{T}(p,
  \theta')\) by post-composition. That is, suppose given open random dynamical
  systems and morphisms over \((p, \theta)\) as in the diagram of Proposition
  \ref{prop:poly-rdyn}. Then \(\RDynT{T}(p, \phi)\) returns the following
  diagram:
  \[\begin{tikzcd}
	X &&&& {\Sum_{x:X} p[\vartheta^o(t, x)]} &&&& X \\
	\\
	&& {\Omega'} &&&& {\Omega'} \\
	\\
	Y &&&& {\Sum_{y:Y} p[\psi^o(t, y)]} &&&& Y
	\arrow["{\theta'(t)}"', from=3-3, to=3-7]
	\arrow["{\vartheta^o(t)^\ast \sigma}", from=1-1, to=1-5]
	\arrow["{\vartheta^u(t)}", from=1-5, to=1-9]
	\arrow["{\psi^o(t)^\ast \sigma}"', from=5-1, to=5-5]
	\arrow["{\psi^u(t)}"', from=5-5, to=5-9]
	\arrow["f"', from=1-1, to=5-1]
	\arrow["f", from=1-9, to=5-9]
	\arrow["{\phi\circ\pi_Y}"', from=5-9, to=3-7]
	\arrow["{\phi\circ\pi_X}", from=1-9, to=3-7]
	\arrow["{\phi\circ\pi_Y}", from=5-1, to=3-3]
	\arrow["{\phi\circ\pi_X}"', from=1-1, to=3-3]
  \end{tikzcd}\]
  That is, \(\RDynT{T}(p, \phi)(\vartheta) := (\phi\circ\pi_X, \vartheta^o,
  \vartheta^u)\) and \(\RDynT{T}(p, \phi)(f)\) is given by the same underlying
  map \(f : X \to Y\) on state spaces.

  \begin{proof}
    We need to check: the naturality condition of Definition \ref{def:poly-rdyn}
    for both \(\RDynT{T}(\varphi, \theta)(\vartheta)\) and \(\RDynT{T}(p,
    \phi)(\vartheta)\); functoriality of \(\RDynT{T}(\varphi, \theta)\) and
    \(\RDynT{T}(p, \phi)\); and (pseudo)functoriality of \(\RDynT{T}\) with
    respect to both morphisms of polynomials and of closed measure-preserving
    systems.

    We begin by checking that the conditions of Definition \ref{def:poly-rdyn}
    and Proposition \ref{prop:poly-dyn-idx} are satisfied by the objects
    \(\RDynT{T}(\varphi, \theta)(\pi_X, \vartheta^o, \vartheta^u) : \RDynT{T}(q,
    \theta)\) and morphisms \(\RDynT{T}(\varphi, \theta)(f) : \RDynT{T}(\varphi,
    \theta)(\pi_X, \vartheta^o, \vartheta^u) \to \RDynT{T}(\varphi,
    \theta)(\pi_Y, \psi^o, \psi^u)\) in the image of \(\RDynT{T}(\varphi,
    \theta)\). We proceed similarly to the proof of Proposition
    \ref{prop:poly-dyn-idx}. Therefore, given a section \(\tau : q(1) \to
    \Sum_{j:q(1)} q[j]\) of q, we need to check that the closure
    \(\RDynT{T}(\varphi, \theta)(\vartheta)^\tau\) forms a closed random
    dynamical system in \(\Cat{Cat}(\deloop{\Tt}, \cat{E})/\theta\). That is to
    say, for all \(t : \Tt\) and sections \(\tau\), we need to check that the
    following naturality square commutes:
    \[\begin{tikzcd}
	X && {\Sum_{x:X} q[\varphi_1 \circ \vartheta^o(t,x)]} && {\Sum_{x:X} p[\vartheta^o(t,x)]} && X \\
	\\
	\Omega &&&&&& \Omega
	\arrow["{\vartheta^o(t)^\ast \tau}", from=1-1, to=1-3]
	\arrow["{\vartheta^o(t)^\ast \varphi^\#}", from=1-3, to=1-5]
	\arrow["{\vartheta^u}", from=1-5, to=1-7]
	\arrow["{\pi_X}"', from=1-1, to=3-1]
	\arrow["{\theta(t)}"', from=3-1, to=3-7]
	\arrow["{\pi_X}", from=1-7, to=3-7]
    \end{tikzcd}\]
    As before, we find that \(\varphi^\# \circ \varphi_1^\ast \tau\) is a
    section of \(p\), so that commutativity of the diagram above reduces to
    commutativity of the diagram in Definition \ref{def:poly-rdyn}. Similarly,
    given a morphism \(f : (\pi_X, \vartheta^o, \vartheta^u) \to (\pi_Y, \psi^o,
    \psi^u)\), we need to check that the diagram in Proposition
    \ref{prop:poly-rdyn} induced for \(\RDynT{T}(\varphi, \theta)(f)\) commutes
    for all times \(t : \Tt\) and sections \(\tau\) of \(q\). But as in the
    proof of Proposition \ref{prop:poly-dyn-idx}, given such a section \(\tau\),
    the diagram for \(\RDynT{T}(\varphi, \theta)(f)\) reduces to that for \(f\)
    and the section \(\varphi^\# \circ \varphi_1^\ast \tau\) of \(p\), which
    commutes \textit{ex hypothesi}; and functoriality of \(\RDynT{T}(\varphi,
    \theta)\) follows immediately.

    Next, we check that the conditions of Definition \ref{def:poly-rdyn} and
    Proposition \ref{prop:poly-dyn-idx} are satisfied in the image of
    \(\RDynT{T}(p, \phi)\). It is clear by the definition of the action of
    \(\RDynT{T}(p, \phi)\) that the condition that the diagram in Proposition
    \ref{prop:poly-dyn-idx} commutes is satisfied, from which it follows by
    pasting that \(\RDynT{T}(p, \phi)\) is functorial. We therefore just have to
    check the induced diagram in Definition \ref{def:poly-rdyn}
    commutes. Consider the following diagram:
    \[\begin{tikzcd}
	X && {\Sum_{x:X} p[\vartheta^o(t, x)]} && X \\
	\\
	\Omega &&&& \Omega \\
	\\
	{\Omega'} &&&& {\Omega'}
	\arrow["{\pi_X}"', from=1-1, to=3-1]
	\arrow["{\vartheta^o(t)^\ast \sigma}", from=1-1, to=1-3]
	\arrow["{\vartheta^u(t)}", from=1-3, to=1-5]
	\arrow["{\pi_X}", from=1-5, to=3-5]
	\arrow["{\theta(t)}", from=3-1, to=3-5]
	\arrow["\phi"', from=3-1, to=5-1]
	\arrow["\phi", from=3-5, to=5-5]
	\arrow["{\theta'(t)}", from=5-1, to=5-5]
    \end{tikzcd}\]
    The top square commutes \textit{ex hypothesi}, the bottom square commutes by
    the definition of morphism of closed measure-preserving dynamical systems
    (Proposition \ref{prop:metric-sys}), and the outer square is the induced
    diagram we need to check, which therefore commutes by the pasting of
    commuting squares.

    Finally, we check that \(\RDynT{T}\) is functorial with respect to morphisms
    of polynomials and morphisms of closed measure-preserving dynamical systems.
    As in the proof of Proposition \ref{prop:poly-dyn-idx}, these reduce to
    checking that pullback and composition are functorial, which we again leave
    to the dedicated reader.
  \end{proof}
\end{prop}

In applications, it is often desirable to connect together systems with
different noise sources, which means collecting together the indexing over
closed measure-preserving systems into a single total category. Formally, this
amounts to constructing an opfibration using the Grothendieck construction.

\begin{prop} \label{prop:poly-rdyn-fib}
  The indexing of \(\RDynT{T}\) by closed measure-preserving systems generates,
  for each polynomial \(p : \Poly{E}\), an opfibration over
  \(\Cat{Cat}(\deloop{\Tt}, \cat{E})_{\Pa}\), denoted \(\int\RDynT{T}(p)\);
  retaining the indexing by polynomials makes its type \(\int\RDynT{T} :
  \Poly{E} \to \Cat{Fib}\big(\Cat{Cat}(\deloop{\Tt}, \cat{E})_{\Pa}\op\big)\).

  Explicitly, an object of \(\int\RDynT{T}(p)\) is a pair \((\vartheta,
  \theta)\) where \(\theta := (\theta^u, \beta)\) is an object of
  \(\Cat{Cat}(\deloop{\Tt}, \cat{E})_{\Pa}\) (\textit{i.e.}, a closed
  measure-preserving dynamical system) and \(\vartheta := (\pi, \vartheta^o,
  \vartheta^u)\) is an object of \(\RDynT{T}(p, \theta)\) (\textit{i.e.}, an
  open random dynamical system over \(\theta\) on the interface \(p\)).  A
  morphism \((\vartheta, \theta) \to (\vartheta', \theta')\) consists in a pair
  \((f, \phi)\), where \(\phi : \theta \to \theta'\) is a morphism of closed
  measure-preserving systems and \(f : \RDynT{T}(p, \phi)(\vartheta) \to
  \vartheta'\) is a morphism in \(\RDynT{T}(p, \theta')\) of open random
  dynamical systems.  The identity morphism on an object of \(\int\RDynT{T}\) is
  given by the corresponding pair of identities.  Given \((f, \phi) :
  (\vartheta, \theta) \to (\vartheta', \theta')\) and \((f', \phi') :
  (\vartheta', \theta') \to (\vartheta'', \theta'')\), their composite is given
  by the pair \(\big(f' \circ \RDynT{T}(p, \phi')(f), \phi' \circ \phi\big)\).

  \begin{proof}
    We start with a doubly-indexed category \(\RDynT{T} : \Poly{E} \times
    \Cat{Cat}(\deloop{\Tt}, \cat{E})_{\Pa} \to \Cat{Cat}\).  By the Cartesian
    closure of \(\Cat{Cat}\), this induces an indexed category of indexed
    categories \(\Poly{E} \to \Cat{Cat}\big(\Cat{Cat}(\deloop{\Tt},
    \cat{E})_{\Pa}, \Cat{Cat}\big)\).  Applying the covariant Grothendieck
    construction to each indexed category in the codomain generates an indexed
    opfibration \(\int\RDynT{T} : \Poly{E} \to
    \Cat{Fib}\big(\Cat{Cat}(\deloop{\Tt}, \cat{E})_{\Pa}\op\big)\). Unpacking
    the structure generated by this indexed Grothendieck construction gives the
    explicit form presented in the proposition.
  \end{proof}
\end{prop}

\begin{prop}
  Suppose \((\theta, \beta)\) is a closed measure-preserving dynamical
  system. As in the deterministic case (Proposition
  \ref{prop:closed-sys-in-dyn-y}), the category \(\RDynT{T}(y, \theta)\) of open
  random dynamical systems on the trivial interface \(y\) with base system
  \(\theta\) is equivalent to the category of closed random dynamical systems
  \(\Cat{Cat}(\deloop{\Tt}, \cat{E})/\theta\).
  \begin{proof}
    The proof is directly analogous to the proof of Proposition
    \ref{prop:closed-sys-in-dyn-y}, except that we now need to check that the
    flow condition induced by Definition \ref{def:poly-rdyn} in the case of the
    trivial interface \(y\) corresponds to the definition of a bundle in
    \(\Cat{Cat}(\deloop{\Tt}, \cat{E})/\theta\), and that the morphisms of
    systems coincide similarly. Both of these verifications are straightforward
    and we leave them to the reader.
  \end{proof}
\end{prop}

To be satisfactorily `open', we should want the noise sources themselves to be
open systems, preserving measure in an appropriately generalized sense.
Fortunately, we can follow the theme of the developments above to define an
indexed category of open measure-preserving systems.

\begin{defn}
  Let \(p : \Poly{E}\) be a polynomial and let \((S, \nu) : 1/\Kl(\Pa)\) be a
  probability space, both in \(\cat{E}\). An \textbf{open measure-preserving
    dynamical system on the interface} \(p\) \textbf{with state space} \((S,
  \nu)\) \textbf{and time} \(\Tt\) consists in a pair of morphisms \(\vartheta^o
  : \Tt \times S \to p(1)\) and \(\vartheta^u : \Sum_{t : \Tt} \Sum_{s : S}
  p[\vartheta^o(t, s)] \to S\), such that, for any global section \(\sigma :
  p(1) \to \Sum_{i:p(1)} p[i]\) of \(p\), the maps \(\vartheta^\sigma : \Tt
  \times S \to S\) given by
  \[
  \Sum_{t:\Tt} S \xto{\vartheta^o(-)^\ast \sigma} \Sum_{t:\Tt} \Sum_{s:S} p[\vartheta^o(-, s)] \xto{\vartheta^u} S
  \]
  constitute a closed measure-preserving dynamical system---again called the
  \textit{closure} of \(\vartheta\)---with respect to \(\nu\), \textit{i.e.}, an
  object in \(\Cat{Cat}(\deloop{\Tt}, \cat{E})_{\Pa}\); see Definition
  \ref{def:metric-sys} for the explicit conditions. We collect the data of such
  an open dynamical system into a tuple \(\vartheta = (S, \nu, \vartheta^o,
  \vartheta^u)\).
\end{defn}

\begin{prop}
  Open measure-preserving dynamical systems over \(p\) with time \(\Tt\) form
  the objects of a category, denoted \(\DynT{T}(p)_{\Pa}\).  Its morphisms \(f :
  (X, \mu, \vartheta^o, \vartheta^u) \to (Y, \nu, \psi^o, \psi^u)\) are maps
  between the state spaces \(f : X \to Y\) in \(\cat{E}\) such that, for all
  times \(t : \Tt\) and sections \(\sigma : p(1) \to \sum_{i:p(1)} p[i]\) of
  \(p\), the map \(f\) lifts to a morphism \(f^\sigma : \vartheta^\sigma \to
  \psi^\sigma\) in \(\Cat{Cat}(\deloop{\Tt}, \cat{E})_{\Pa}\) between the
  closures. More explicitly, this condition corresponds to an extension of the
  condition on morphisms of plain open dynamical systems (Proposition
  \ref{prop:poly-dyn}) so that they preserve measure (Proposition
  \ref{prop:metric-sys}).
  \begin{proof}
    The proof amounts to the straightforward proof of Proposition
    \ref{prop:poly-dyn}, extended so that morphisms and their composition
    preserve measure. Just as checking compositionality in the plain case is
    straightforward, so is checking that measure is preserved: the relevant
    composite triangles must commute by pasting.
  \end{proof}
\end{prop}

\begin{prop} \label{prop:poly-metric-idx}
  The categories \(\DynT{T}(p)_{\Pa}\) collect into a polynomially-indexed
  category, \(\DynT{T}_{\Pa} : \Poly{E} \to \Cat{Cat}\). The action of
  \(\DynT{T}_{\Pa}\) on morphisms \(\varphi : p \to q\) of polynomials is
  defined as for plain open dynamical systems in Proposition
  \ref{prop:poly-dyn-idx}.
  \begin{proof}
    The proof proceeds as for the proof of Proposition
    \ref{prop:poly-dyn-idx}. We also need to check that the structures in the
    image of \(\DynT{T}_{\Pa}(\varphi)\) for each polynomial morphism \(\varphi
    : p \to q\) satisfy the relevant measure-preservation property. But as in
    the other parts of the proof, this property follows from the facts that the
    structures in \(\DynT{T}(p)_{\Pa}\) satisfy the property \textit{ex
      hypothesi} for any section of \(p\), and that any section \(\tau\) of
    \(q\) can be pulled back along \(\varphi\) to a section of \(p\)
    accordingly.
  \end{proof}
\end{prop}

\begin{prop} \label{prop:metric-in-open-metric}
  Just as the category \(\Cat{Cat}(\deloop{\Tt}, \cat{E})\) of general closed
  dynamical systems is equivalent to \(\DynT{T}(y)\) (Proposition
  \ref{prop:closed-sys-in-dyn-y}), there is an analogous equivalence in the
  measure-preserving case; that is, \(\Cat{Cat}(\deloop{\Tt}, \cat{E})_{\Pa}
  \cong \DynT{T}(y)_{\Pa}\).
  \begin{proof}
    Given Proposition \ref{prop:closed-sys-in-dyn-y}, we only need to check that
    the measure-preservation condition coincides for objects and morphisms; but
    this is immediate from the definitions.
  \end{proof}
\end{prop}

\begin{prop}
  There is a forgetful indexed functor \(U : \DynT{T}_{\Pa} \to \DynT{T}\) which
  simply forgets the probability space structures.
\end{prop}

Random dynamical systems form a subcategory of bundles of dynamical systems, and
so, to define a general category of ``fully open'' random dynamical systems, it
makes sense to start by defining open bundles of dynamical systems.

\begin{defn} \label{defn:bdyn-pbth}
  Let \(p, b : \Poly{E}\) be polynomials in \(\cat{E}\), and let \(\theta :=
  (\theta(\ast), \theta^o, \theta^u) : \DynT{T}(b)\) be an open dynamical system
  over \(b\).  An \textbf{open bundle dynamical system} over \((p, b, \theta)\)
  is a pair \((\pi_{\vartheta\theta}, \vartheta)\) where \(\vartheta :=
  (\vartheta(\ast), \vartheta^o, \vartheta^u) : \DynT{T}(p)\) is an open
  dynamical system over \(p\) and \(\pi_{\vartheta\theta} : \vartheta(\ast) \to
  \theta(\ast)\) is a bundle in \(\cat{E}\), such that, for all time \(t : \Tt\)
  and sections \(\sigma\) of \(p\) and \(\varsigma\) of \(b\), the following
  diagrams commute, thereby inducing a bundle of closed dynamical systems
  \(\pi^{\sigma\varsigma}_{\vartheta\theta} : \vartheta^\sigma \to
  \theta^\varsigma\) in \(\Cat{Cat}(\deloop{\Tt}, \cat{E})\):
  \[\begin{tikzcd}
	\vartheta(\ast) && {\Sum_{w:\vartheta(\ast)} p[\vartheta^o(t, w)]} && {\vartheta(\ast)} \\
	\\
	\theta(\ast) && {\Sum_{x:\theta(\ast)} b[\theta^o(t,x)]} && \theta(\ast)
	\arrow["{\pi_{\vartheta\theta}}", from=1-1, to=3-1]
	\arrow["{\pi_{\vartheta\theta}}", from=1-5, to=3-5]
	\arrow["{\vartheta^o(t)^\ast \sigma}", from=1-1, to=1-3]
	\arrow["{\vartheta^u(t)}", from=1-3, to=1-5]
	\arrow["{\theta^o(t)^\ast \varsigma}", from=3-1, to=3-3]
	\arrow["{\theta^u(t)}", from=3-3, to=3-5]
  \end{tikzcd}\]
\end{defn}

\begin{prop} \label{prop:bdyn-pbth-cat}
  Let \(p, b : \Poly{E}\) be polynomials in \(\cat{E}\), and let \(\theta :=
  (\theta(\ast), \theta^o, \theta^u) : \DynT{T}(b)\) be an open dynamical system
  over \(b\).  Open bundle dynamical systems over \((p, b, \theta)\) form the
  objects of a category \(\BDynT{T}(p, b, \theta)\). Morphisms \(f :
  (\pi_{\vartheta\theta}, \vartheta) \to (\pi_{\varrho\theta}, \varrho)\) are
  maps \(f : \vartheta(\ast) \to \varrho(\ast)\) in \(\cat{E}\) making the
  following diagram commute for all times \(t : \Tt\) and sections \(\sigma\) of
  \(p\) and \(\varsigma\) of \(b\):
  \[\begin{tikzcd}
	{\vartheta(\ast)} &&&& {\Sum_{w:\vartheta(\ast)} p[\vartheta^o(t, w)]} &&&& {\vartheta(\ast)} \\
	\\
	&& {\theta(\ast)} && {\Sum_{x:\theta(\ast)} b[\theta^o(t,x)]} && {\theta(\ast)} \\
	\\
	{\varrho(\ast)} &&&& {\Sum_{y:\varrho(\ast)} p[\varrho^o(t, y)]} &&&& {\varrho(\ast)}
	\arrow["{\vartheta^o(t)^\ast \sigma}", from=1-1, to=1-5]
	\arrow["{\vartheta^u(t)}", from=1-5, to=1-9]
	\arrow["{\theta^o(t)^\ast \varsigma}", from=3-3, to=3-5]
	\arrow["{\theta^u(t)}", from=3-5, to=3-7]
	\arrow["{\pi_{\vartheta\theta}}"', from=1-1, to=3-3]
	\arrow["{\pi_{\varrho\theta}}", from=5-1, to=3-3]
	\arrow["{\pi_{\vartheta\theta}}", from=1-9, to=3-7]
	\arrow["{\pi_{\varrho\theta}}"', from=5-9, to=3-7]
	\arrow["{\varrho^o(t)^\ast \sigma}"', from=5-1, to=5-5]
	\arrow["{\varrho^u(t)}"', from=5-5, to=5-9]
	\arrow["{f}"', from=1-1, to=5-1]
	\arrow["{f}", from=1-9, to=5-9]
  \end{tikzcd}\]
  That is, \(f\) is a map on the state spaces that induces a morphism
  \((\pi_{\vartheta\theta}, \vartheta^\sigma) \to (\pi_{\varrho\theta},
  \varrho^\sigma)\) in \(\Cat{Cat}(\deloop{\Tt}, \cat{E})/\theta^\varsigma\) of
  bundles of the closures. Identity morphisms are the corresponding identity
  maps, and composition is by pasting.
\end{prop}

\begin{prop} \label{prop:bdyn-bth-idx}
  Varying the polynomials \(p\) in \(\BDynT{T}(p, b, \theta)\) induces a
  polynomially indexed category \(\BDynT{T}({-}, b, \theta) : \Poly{E} \to
  \Cat{Cat}\). On polynomials \(p\), it returns the categories \(\BDynT{T}(p, b,
  \theta)\) of Proposition \ref{prop:bdyn-pbth-cat}. On morphisms \(\varphi : p
  \to q\) of polynomials, define the functors \(\BDynT{T}(\varphi, b, \theta) :
  \BDynT{T}(p, b, \theta) \to \BDynT{T}(q, b, \theta)\) as in Propositions
  \ref{prop:poly-dyn-idx} and \ref{prop:poly-rdyn-idx}. That is, suppose
  \((\pi_{\vartheta\theta}, \vartheta) : \BDynT{T}(p, b, \theta)\) is object
  (open bundle dynamical system) in \(\BDynT{T}(p, b, \theta)\), where
  \(\vartheta := (\vartheta(\ast), \vartheta^o, \vartheta^u)\). Then its image
  \(\BDynT{T}(\varphi, b, \theta)(\pi_{\vartheta\theta}, \vartheta)\) is defined
  as the pair \((\pi_{\vartheta\theta}, \varphi\vartheta)\), where
  \(\varphi\vartheta := (\vartheta(\ast), \phi_1 \circ \vartheta^o, \vartheta^u
  \circ {\vartheta^o}^\ast \varphi^\#)\). On morphisms \(f :
  (\pi_{\vartheta\theta}, \vartheta) \to (\pi_{\varrho\theta}, \varrho)\),
  \(\BDynT{T}(\varphi, b, \theta)(f)\) is again given by the same underlying map
  \(f : \vartheta(\ast) \to \varrho(\ast)\) of state spaces.
  \begin{proof}
    The proof amounts to the proof for Proposition \ref{prop:poly-rdyn-idx} that
    \(\RDynT{T}(\varphi, \theta)\) constitutes an indexed category, except that
    the closed base dynamical system \(\theta\) of that Proposition is here
    replaced, for any section \(\varsigma\) of \(b\), by the closure
    \(\theta^\varsigma\) by \(\varsigma\) of the open dynamical system \(\theta
    : \DynT{T}(b)\) of the present Proposition. The proof goes through
    accordingly, since the relevant diagrams are guaranteed to commute for any
    such \(\varsigma\) by the conditions in Definition \ref{defn:bdyn-pbth} and
    Proposition \ref{prop:bdyn-pbth-cat}.
  \end{proof}
\end{prop}

\begin{prop} \label{prop:bdyn-b-idx}
  Letting the base system \(\theta\) also vary induces a doubly-indexed category
  \(\BDynT{T}({-}, b, {=}) : \Poly{E} \times \DynT{T}(b) \to \Cat{Cat}\). Given
  a polynomial \(p : \Poly{E}\) and morphism \(\phi : \theta \to \rho\) in
  \(\DynT{T}(b)\), the functor \(\BDynT{T}(p, b, \phi) : \BDynT{T}(p, b, \theta)
  \to \BDynT{T}(p, b, \rho)\) is defined by post-composition, as in Proposition
  \ref{prop:poly-rdyn-idx} for the action of \(\RDynT{T}\) on morphisms of the
  base systems there. More explicitly, such a morphism \(\phi\) corresponds to a
  map \(\phi : \theta(\ast) \to \rho(\ast)\) of state spaces in
  \(\cat{E}\). Given an object \((\pi_{\vartheta\theta}, \vartheta)\) of
  \(\BDynT{T}(p, b, \theta)\), we define \(\BDynT{T}(p, b,
  \phi)(\pi_{\vartheta\theta}, \vartheta) := (\phi \circ \pi_{\vartheta\theta},
  \vartheta)\). Given a morphism \(f : (\pi_{\vartheta\theta}, \vartheta) \to
  (\pi_{\varrho\theta}, \varrho)\) in \(\BDynT{T}(p, b, \theta)\), its image
  \(\BDynT{T}(p, b, \phi)(f) : (\phi\circ\pi_{\vartheta\theta}, \vartheta) \to
  (\phi\circ\pi_{\varrho\theta}, \varrho)\) is given by the same underlying map
  \(f : \vartheta(\ast) \to \varrho(\ast)\) of state spaces.
  \begin{proof}
    As for Proposition \ref{prop:bdyn-bth-idx}, the proof here amounts to the
    proof for Proposition \ref{prop:poly-rdyn-idx} that \(\RDynT{T}(p, \phi)\)
    constitutes an indexed category, except again the closed systems are
    replaced by (the appropriate closures of) open ones, and the
    measure-preserving structure is forgotten.
  \end{proof}
\end{prop}

\begin{prop} \label{prop:bdyn-b-fib}
  There is an indexed opfibration \(\int\BDynT{T}({-}, b) : \Poly{E} \to
  \Cat{Fib}\big(\DynT{T}(b)\big)\) generated from \(\BDynT{T}({-}, b, {=})\) by
  the Grothendieck construction.

  Explicitly, given a polynomial \(p : \Poly{E}\), the objects of
  \(\int\BDynT{T}(p, b)\) are triples \((\pi_{\vartheta\theta}, \vartheta,
  \theta)\), where \(\theta : \DynT{T}(b)\) is an open dynamical system over
  \(b\) and \((\pi_{\vartheta\theta}, \vartheta) : \BDynT{T}(p, b, \theta)\) is
  an open bundle dynamical system over \((p, b, \theta)\). Morphisms \(f :
  (\pi_{\vartheta\theta}, \vartheta, \theta) \to (\pi_{\varrho\rho}, \varrho,
  \rho)\) are pairs \((f_p, f_b)\) of a morphism \(f_b : \theta \to \rho\) in
  \(\DynT{T}(b)\) and a morphism \(f_p : (f_b \circ \pi_{\vartheta\theta},
  \vartheta) \to (\pi_{\varrho\rho}, \varrho)\) in \(\BDynT{T}(p, b, \rho)\)
  making the following diagram commute for all sections \(\sigma\) of \(p\) and
  \(\varsigma\) of \(b\):
  \[\begin{tikzcd}
	{\vartheta(\ast)} && {\Sum_{w:\vartheta(\ast)} p[\vartheta^o(t, w)]} && {\vartheta(\ast)} \\
	\\
	& {\theta(\ast)} && {\Sum_{x:\theta(\ast)} b[\theta^o(t, x)]} && {\theta(\ast)} \\
	\\
	{\varrho(\ast)} && {\Sum_{y:\varrho(\ast)} p[\varrho^o(t,y)]} && {\varrho(\ast)} \\
	\\
	& {\rho(\ast)} && {\Sum_{z:\rho(\ast)} b[\rho^o(t, z)]} && {\rho(\ast)}
	\arrow["{\pi_{\vartheta\theta}}", from=1-1, to=3-2]
	\arrow["{\pi_{\vartheta\theta}}", from=1-5, to=3-6]
	\arrow["{\vartheta^o(t)^\ast \sigma}", from=1-1, to=1-3]
	\arrow["{\theta^o(t)^\ast \varsigma}", from=3-2, to=3-4]
	\arrow["{\vartheta^u(t)}", from=1-3, to=1-5]
	\arrow["{\theta^u(t)}"{pos=0.3}, from=3-4, to=3-6]
	\arrow["{\varrho^o(t)^\ast \sigma}"{pos=0.7}, from=5-1, to=5-3]
	\arrow["{\varrho^u(t)}", from=5-3, to=5-5]
	\arrow["{\rho^o(t)^\ast \varsigma}", from=7-2, to=7-4]
	\arrow["{\rho^u(t)}", from=7-4, to=7-6]
	\arrow["{\pi_{\varrho\rho}}"', from=5-1, to=7-2]
	\arrow["{\pi_{\varrho\rho}}"', from=5-5, to=7-6]
	\arrow["{f_p}"'{pos=0.4}, curve={height=6pt}, from=1-1, to=5-1]
	\arrow["{f_b}"'{pos=0.3}, curve={height=6pt}, from=3-2, to=7-2]
	\arrow["{f_p}"{pos=0.7}, curve={height=-6pt}, from=1-5, to=5-5]
	\arrow["{f_b}"{pos=0.6}, curve={height=-6pt}, from=3-6, to=7-6]
  \end{tikzcd}\]
  Identity morphisms are the pairs of the corresponding identities, and
  composition is again by pasting.
  \begin{proof}
    Compare Proposition \ref{prop:poly-rdyn-fib}.
  \end{proof}
\end{prop}

\begin{prop} \label{prop:bdyn-idx}
  Varying the base polynomial \(b\) extends \(\BDynT{T}({-}, b, {=})\) to a
  triply indexed category, \(\BDynT{T} : \Poly{E} \times \sum_{b:\Poly{E}}
  \DynT{T}(b) \to \Cat{Cat}\) and \(\int\BDynT{T}({-}, b)\) to a doubly indexed
  fibration, \(\int\BDynT{T} : \Poly{E} \to \prod_{b:\Poly{E}}
  \Cat{Fib}\big(\DynT{T}(b)\big)\); these are equivalent by the Grothendieck
  construction.

  Let \(p, b, c : \Poly{E}\) be polynomials, and let \(\chi : b \to c\) be a
  morphism accordingly. Let \(\theta : \DynT{T}(b)\) range over open dynamical
  systems over \(b\). We define the (dependent) functor \(\BDynT{T}(p, \chi,
  \theta) : \BDynT{T}(p, b, \theta) \allowbreak\to \BDynT{T}(p, c,
  \DynT{T}(\chi)(\theta))\) as follows.  This functor is equivalent to its image
  under the Grothendieck construction, \(\int\BDynT{T}(p, \chi) :
  \int\BDynT{T}(p, b) \to \int\BDynT{T}(p, c)\), which is easier to
  describe. Therefore, let \((\pi_{\vartheta\theta}, \vartheta, \theta)\) be an
  object of \(\int\BDynT{T}(p, b)\). Its image under \(\int\BDynT{T}(p, \chi)\)
  is the object \((\pi_{\vartheta\theta}, \vartheta, \DynT{T}(\chi)(\theta)) :
  \int\BDynT{T}(p, c)\). Given a morphism \((f_p, f_b) : (\pi_{\vartheta\theta},
  \vartheta, \theta) \to (\pi_{\varrho\rho}, \varrho, \rho)\) in
  \(\int\BDynT{T}(p, b)\), its image \(\int\BDynT{T}(p, \chi)(f) :
  (\pi_{\vartheta\theta}, \vartheta, \DynT{T}(\chi)(\theta)) \to
  (\pi_{\varrho\rho}, \varrho, \DynT{T}(\chi)(\rho))\) is given by the same maps
  \(f_p : \vartheta(\ast) \to \varrho(\ast)\) and \(f_b : \theta(\ast) \to
  \rho(\ast)\) of state spaces.
  \begin{proof}
    We just need to check that the diagrams of Definition \ref{defn:bdyn-pbth}
    and Proposition \ref{prop:bdyn-b-fib} induced in the image of each
    \(\int\BDynT{T}(p, \chi)\) commute, and that \(\int\BDynT{T}(p, {-})\) is
    functorial with respect to morphisms of polynomials.  To check the first
    diagram, we note that the following commutes for the usual reason, that
    \(\chi^\# \circ \chi_1^\ast \tau\) is a section of \(b\):
    \[\begin{tikzcd}
	{\vartheta(\ast)} &&& {\Sum_{w:\vartheta(\ast)} p[\vartheta^o(t, w)]} &&& {\vartheta(\ast)} \\
	\\
	{\theta(\ast)} && {\Sum_{x:\theta(\ast)} c[\chi_1 \circ \theta^o(t, x)]} && {\Sum_{x:\theta(\ast)} b[\theta^o(t,x)]} && {\theta(\ast)}
	\arrow["{\vartheta^o(t)^\ast \sigma}", from=1-1, to=1-4]
	\arrow["{\vartheta^u(t)}", from=1-4, to=1-7]
	\arrow["{\pi_{\vartheta\theta}}", from=1-7, to=3-7]
	\arrow["{\pi_{\vartheta\theta}}"', from=1-1, to=3-1]
	\arrow["{\theta^o(t)^\ast \chi_1^\ast \tau}", from=3-1, to=3-3]
	\arrow["{\theta^o(t)^\ast \chi^\#}", from=3-3, to=3-5]
	\arrow["{\theta^u(t)}", from=3-5, to=3-7]
    \end{tikzcd}\]
    To check the second diagram, we observe that the same property makes the
    following commute, given a morphism \((f_p, f_b) : (\pi_{\vartheta\theta},
    \vartheta, \theta) \to (\pi_{\varrho\rho}, \varrho, \rho)\) in
    \(\int\BDynT{T}(p, b)\):
    \[\begin{tikzcd}
	{\vartheta(\ast)} &&& {\Sum_{w:\vartheta(\ast)} p[\vartheta^o(t,w)]} &&& {\vartheta(\ast)} \\
	\\
	& {\theta(\ast)} && {\Sum_{x:\theta(\ast)} c[\chi_1\circ\theta^o(t,x)]} && {\Sum_{x:\theta(\ast)}b[\theta^o(t,x)]} && {\theta(\ast)} \\
	\\
	{\varrho(\ast)} &&& {\Sum_{y:\varrho(\ast)} p[\varrho^o(t,y)]} &&& {\varrho(\ast)} \\
	\\
	& {\rho(\ast)} && {\Sum_{z:\rho(\ast)} c[\chi_1\circ\rho^o(t,z)]} && {\Sum_{z:\rho(\ast)} b[\rho^o(t,z)]} && {\rho(\ast)}
	\arrow["{\vartheta^o(t)^\ast \sigma}", from=1-1, to=1-4]
	\arrow["{\vartheta^u(t)}", from=1-4, to=1-7]
	\arrow["{\pi_{\vartheta\theta}}", from=1-1, to=3-2]
	\arrow["{\theta^o(t)^\ast \chi_1^\ast \tau}", from=3-2, to=3-4]
	\arrow["{\theta^o(t)^\ast \chi^\#}", from=3-4, to=3-6]
	\arrow["{\theta^u(t)}"{pos=0.3}, from=3-6, to=3-8]
	\arrow["{\pi_{\varrho\rho}}", from=1-7, to=3-8]
	\arrow["{f_p}"'{pos=0.4}, curve={height=6pt}, from=1-1, to=5-1]
	\arrow["{f_b}"'{pos=0.3}, curve={height=6pt}, from=3-2, to=7-2]
	\arrow["{\pi_{\vartheta\theta}}"'{pos=0.4}, from=5-1, to=7-2]
	\arrow["{\pi_{\varrho\rho}}"'{pos=0.4}, from=5-7, to=7-8]
	\arrow["{f_p}"{pos=0.7}, curve={height=-6pt}, from=1-7, to=5-7]
	\arrow["{f_b}"{pos=0.6}, curve={height=-6pt}, from=3-8, to=7-8]
	\arrow["{\varrho^o(t)^\ast \sigma}"{pos=0.6}, from=5-1, to=5-4]
	\arrow["{\varrho^u(t)}", from=5-4, to=5-7]
	\arrow["{\rho^o(t)^\ast \chi_1^\ast \tau}", from=7-2, to=7-4]
	\arrow["{\rho^o(t)^\ast \chi^\#}", from=7-4, to=7-6]
	\arrow["{\rho^u(t)}", from=7-6, to=7-8]
    \end{tikzcd}\]
    Finally, we note that, as in the proof of Proposition
    \ref{prop:poly-dyn-idx}, functoriality on morphisms \(b \to c \to d\) of
    polynomials follows from the functoriality of pullback and composition.
  \end{proof}
\end{prop}

\begin{prop} \label{prop:rbdyn-idx}
  By restricting \(\BDynT{T}\) and \(\int\BDynT{T}\) to those systems which
  preserve measure in the base, we obtain the indexed categories \(\RBDynT{T} :
  \Poly{E} \times \sum_{b:\Poly{E}} \DynT{T}(b)_{\Pa} \to \Cat{Cat}\) and
  \(\int\RBDynT{T} : \Poly{E} \to \prod_{b:\Poly{E}}
  \Cat{Fib}\big(\DynT{T}(b)_{\Pa}\big)\) of \textbf{open random bundle dynamical
    systems}, or alternatively, \textbf{fully open random dynamical systems}.
  \begin{proof}
    The only extra check required is that the measure-preservation property is
    retained by the action of \(\RBDynT{T}\) on morphisms of polynomials \(p \to
    q\) above and below \(b \to c\) and on morphisms of open measure-preserving
    systems in the base.  On the latter, the result is immediate, since such
    morphisms are defined to preserve measure.  On morphisms of base polynomials
    \(b\), the result follows from the fact that \(\DynT{T}_{\Pa}\) forms an
    indexed category (Proposition \ref{prop:poly-metric-idx}).  Finally,
    morphisms of polynomials \(p\) above \(b\) induce no change in the base
    systems, and so the measure-preservation property is certainly still
    satisfied.
  \end{proof}
\end{prop}

\begin{prop} \label{prop:rdyn-in-rbdyn}
  Open random dynamical systems embed into fully open random dynamical systems
  with trivial base polynomial. That is, \(\RDynT{T} \cong \RBDynT{T}(y) :
  \Poly{E} \times \DynT{T}(y)_{\Pa} \to \Cat{Cat}\).
  \begin{proof}
    This follows immediately from the facts in Proposition
    \ref{prop:metric-in-open-metric} that \(\DynT{T}(y)_{\Pa}\) is equivalent to
    the category of closed measure-preserving systems, and that open random
    dynamical systems are fully open random dynamical systems with closed bases.
  \end{proof}
\end{prop}

When the polynomial \(p\) is nested over \(b\), meaning that \(p : \Poly{E}(b)\)
in the sense of Proposition \ref{prop:poly-nest}, we should expect the dynamics
on the two interfaces to be compatible with this nesting. We can formalize this
with the following structure.

\begin{prop} \label{prop:ndyn-idx}
  There is an indexed category \(\NDynT{T} : \sum_{b : \Poly{E}(y)} \Poly{E}(b)
  \times \DynT{T}(b) \to \Cat{Cat}\) of \textbf{open dynamical systems over
    nested polynomials}, and a corresponding dependently-indexed fibration
  \(\int\NDynT{T} : \prod_{b:\Poly{E}(y)} \Poly{E}(b) \to
  \Cat{Fib}\big(\DynT{T}(b)\big)\). These are defined as \(\BDynT{T}\) and
  \(\int\BDynT{T}\) with the extra condition that the dynamics are compatible
  with the nesting. More explicitly, suppose that \(p\) is a polynomial nested
  over \(b\). Then there are morphisms \(m : \sum_i p[i] \to \sum_k b[k]\) and
  \(n : p(1) \to b(1)\) in \(\cat{E}\), such that \(n \circ p = b \circ m\). An
  open dynamical system over \((m, n) : p \to b\), \textit{i.e.} an object of
  \(\int\NDynT{T}(m,n) : \Cat{Fib}\big(\DynT{T}(b)\big)\), is a triple
  \((\pi_{\vartheta\theta}, \vartheta, \theta)\) where \(\theta : \DynT{T}(b)\)
  is an open dynamical system over \(b\), \(\vartheta : \DynT{T}(p)\) is an open
  dynamical system over \(p\), and \(\pi_{\vartheta\theta} : \vartheta(\ast) \to
  \theta(\ast)\) is a map between the state spaces which lifts uniquely to make
  the front face of the following cube, and therefore the whole cube, commute:
  \[\begin{tikzcd}
	&& {\Sum_{i:p(1)} p[i]} &&& {\Sum_{k:b(1)} b[k]} \\
	\\
	{\Sum_{w:\vartheta(\ast)} p[\vartheta(t, w)]} &&& {\Sum_{x:\theta(\ast)} b[\theta(t,x)]} \\
	& {} & {p(1)} &&& {b(1)} \\
	\\
	{\vartheta(\ast)} &&& {\theta(\ast)}
	\arrow["{\pi_{\vartheta\theta}}", from=6-1, to=6-4]
	\arrow[dashed, from=3-1, to=3-4]
	\arrow[from=3-1, to=6-1]
	\arrow[from=3-4, to=6-4]
	\arrow[from=3-1, to=1-3]
	\arrow[from=3-4, to=1-6]
	\arrow["m", from=1-3, to=1-6]
	\arrow[from=1-3, to=4-3]
	\arrow[from=1-6, to=4-6]
	\arrow["n", from=4-3, to=4-6]
	\arrow["{\vartheta^o(t)}"'{pos=0.6}, from=6-1, to=4-3]
	\arrow["{\theta^o(t)}"', from=6-4, to=4-6]
  \end{tikzcd}\]
  We call the above condition the \textit{nesting condition}; the unlabelled
  edges of this cube are the obvious projections.  As usual, we also require an
  open dynamical system \((\pi_{\vartheta\theta}, \vartheta, \theta)\) to
  satisfy a flow condition, such that for all times \(t : \Tt\) and sections
  \(\sigma\) of \(p\) and \(\varsigma\) of \(b\), the following diagrams
  commute, where the dashed arrow below is the same dashed lift above:
  \[\begin{tikzcd}
	{\vartheta(\ast)} && {\Sum_{w:\vartheta(\ast)} p[\vartheta^o(t,w)]} && {\vartheta(\ast)} \\
	\\
	{\theta(\ast)} && {\Sum_{x:\theta(\ast)} b[\theta^o(t,x)]} && {\theta(\ast)}
	\arrow["{\pi_{\vartheta\theta}}", from=1-1, to=3-1]
	\arrow["{\pi_{\vartheta\theta}}", from=1-5, to=3-5]
	\arrow["{\vartheta^o(t)^\ast \sigma}", from=1-1, to=1-3]
	\arrow["{\sigma^o(t)^\ast \varsigma}"', from=3-1, to=3-3]
	\arrow["{\vartheta^u(t)}", from=1-3, to=1-5]
	\arrow["{\theta^u(t)}"', from=3-3, to=3-5]
	\arrow[dashed, from=1-3, to=3-3]
  \end{tikzcd}\]
  Morphisms in \(\int\NDynT{T}(m, n)\) are as for \(\int\BDynT{T}\) (Proposition
  \ref{prop:bdyn-b-fib}), with the addition of the dashed arrows on the top and
  bottom of the defining cube (in \ref{prop:bdyn-b-fib}). The actions of
  \(\NDynT{T}\) and of \(\int\NDynT{T}\) on morphisms of polynomials, nested
  polynomials, and base systems are defined as for \(\BDynT{T}\); given a
  morphism of polynomials, the dashed lifts are transformed by pullback. Because
  the pasting of two pullback squares is again a pullback square, it is easy to
  check that this also constitutes an indexed category.
\end{prop}

The existence of the dashed lift asserts that the bundle of state spaces
\(\vartheta(\ast) \to \theta(\ast)\) is compatible with the nesting of
polynomials \(p \to b\), in the sense that each section \(\sigma\) of \(p\)
projects onto a compatible section \(\varsigma\) of \(b\).

\begin{prop}
  By restricting \(\NDynT{T}\) and \(\int\NDynT{T}\) to those systems which
  preserve measure in the base, we obtain the indexed categories \(\RNDynT{T} :
  \sum_{b : \Poly{E}(y)} \Poly{E}(b) \times \DynT{T}(b)_{\Pa} \to \Cat{Cat}\)
  and \(\int\RNDynT{T} : \prod_{b:\Poly{E}(y)} \Poly{E}(b) \to
  \Cat{Fib}\big(\DynT{T}(b)_{\Pa}\big)\) of \textbf{random nested dynamical
    systems}.
  \begin{proof}
    The construction and proof are analogous to those of Proposition
    \ref{prop:rbdyn-idx}, with consideration for the nesting condition of
    Proposition \ref{prop:ndyn-idx}.
  \end{proof}
\end{prop}

\begin{prop} \label{prop:rdyn-in-rndyn}
  Open random dynamical systems embed into random nested dynamical systems with
  trivial base polynomial. That is, \(\RDynT{T} \cong \RNDynT{T}(y) :
  \Poly{E}(y) \times \DynT{T}(y)_{\Pa} \to \Cat{Cat}\).
  \begin{proof}
    The result is the same as Proposition \ref{prop:rdyn-in-rbdyn}, with the
    additional nesting condition.  We therefore need to check that there exists
    a unique dashed lift as in Proposition \ref{prop:ndyn-idx}.  Since we have
    \(b = y\) \textit{ex hypothesi}, we have \(\sum_{x:\theta(\ast)}
    b[\theta^o(t, x)] \cong \sum_{x:\theta(\ast)} 1 \cong \theta(\ast)\). Then,
    letting \(\lambda\) denote the projection \(\sum_{w:\vartheta(\ast)}
    p[\vartheta^o(t, w)] \to \vartheta(\ast)\), the dashed lift must be equal to
    \(\pi_{\vartheta\theta} \circ \lambda\).
  \end{proof}
\end{prop}

\begin{ex}
  Let \(p\) be a polynomial, \(M : \cat{E}\) an object, and \(g : M \to p(1)\) a
  map. Then suppose \(\d x_t = f(t, x_t, a_{x_t}) \d t + \sigma(t, x_t) \d W_t\)
  is a stochastic differential equation, with \(f : \rr_+ \times \sum_{x:M}
  p[g(x)] \to TM\). Its solutions \(\chi : \rr_+ \times \Omega \times \sum_{x:M}
  p[g(x)] \to M\) induce an open random dynamical system \((\pi_\Omega : \Omega
  \times M \to M, g, \chi)\) on the interface \(p\) with Wiener base flow
  \((\theta, \gamma)\), following the recipe in Example \ref{ex:brown-sde}.
\end{ex}

\subsection{Internalizing dynamics in Poly}
\label{sec:orgcd5bd98}

\begin{prop} \label{prop:dyn-n-coalg}
  When \(\Tt = \nn\), the category \(\DynT{N}(p)\) of open dynamical systems
  over \(p\) with time \(\nn\) is equivalent to the topos \(p\mdash\Cat{Coalg}\)
  of \(p\)-coalgebras.
  \begin{proof}
    \(p\mdash\Cat{Coalg}\) has as objects pairs \((S, \beta)\) where \(S :
    \cat{E}\) is an object in \(\cat{E}\), \(\beta : S \to p \triangleleft S\)
    is a morphism of polynomials (interpreting \(S\) as the constant copresheaf
    on the set \(S\)), and \(\triangleleft\) denotes the composition monoidal
    product in \(\Poly{E}\) (\textit{i.e.}, composing the corresponding
    copresheaves \(\cat{E} \to \cat{E}\)). A straightforward computation shows
    that, interpreted as an object in \(\cat{E}\), \(p \triangleleft S\)
    corresponds to \(\sum_{i:p(1)} S^{\,p[i]}\).  By the universal property of
    the dependent sum, a morphism \(\beta : S \to \sum_{i:p(1)} S^{\,p[i]}\)
    therefore corresponds bijectively to a pair of maps \(\beta^o : S \to p(1)\)
    and \(\beta^u : \sum_{s:S} p[\beta^o(s)] \to X\).  By Proposition
    \ref{prop:open-transition-map}, such a pair is equivalently a discrete-time
    open dynamical system over \(p\) with state space \(S\): that is, the
    objects of \(p\mdash\Cat{Coalg}\) are in bijection with those of
    \(\DynT{N}(p)\).

    Next, we show that the hom-sets \(p\mdash\Cat{Coalg}\big((S, \beta), (S',
    \beta')\big)\) and \(\DynT{N}(p)\big((S, \beta^o, \beta^u), (S', \beta'^o,
    \beta'^u)\big)\) are in bijection.  A morphism \(f : (S, \beta) \to (S',
    \beta')\) of \(p\)-coalgebras is a morphism \(f : S \to S'\) between the
    state spaces such that \(\beta' \circ f = (p \triangleleft f) \circ
    \beta\). Unpacking this, we find that this means the following diagram in
    \(\cat{E}\) must commute for any section \(\sigma\) of \(p\):
    \[\begin{tikzcd}
	S && {\Sum_{s:S} p[\beta^o(s)]} && {\Sum_{i:p(1)} p[i]} && {p(1)} \\
	&&& {} \\
	&& S && {p(1)} \\
	\\
	&& {S'} && {p(1)} \\
	&&& {} \\
	{S'} && {\Sum_{s':S'} p[\beta'^o(s')]} && {\Sum_{i:p(1)} p[i]} && {p(1)}
	\arrow[Rightarrow, no head, from=3-5, to=5-5]
	\arrow["f"', from=1-1, to=7-1]
	\arrow["f"', from=3-3, to=5-3]
	\arrow[from=1-3, to=1-5]
	\arrow["{\beta^o}", from=3-3, to=3-5]
	\arrow[from=1-5, to=3-5]
	\arrow[from=1-3, to=3-3]
	\arrow["\lrcorner"{anchor=center, pos=0.125}, draw=none, from=1-3, to=2-4]
	\arrow[from=7-3, to=7-5]
	\arrow[from=7-5, to=5-5]
	\arrow["{\beta'^o}", from=5-3, to=5-5]
	\arrow[from=7-3, to=5-3]
	\arrow["\lrcorner"{anchor=center, pos=0.125, rotate=90}, draw=none, from=7-3, to=6-4]
	\arrow["{\beta^u}"', from=1-3, to=1-1]
	\arrow["{\beta'^u}", from=7-3, to=7-1]
	\arrow[Rightarrow, no head, from=3-5, to=1-7]
	\arrow[Rightarrow, no head, from=5-5, to=7-7]
	\arrow["\sigma", from=7-7, to=7-5]
	\arrow["\sigma"', from=1-7, to=1-5]
    \end{tikzcd}\]
    Pulling the arbitrary section \(\sigma\) back along the `output' maps
    \(\beta^o\) and \(\beta'^o\) means that the following commutes:
    \[\begin{tikzcd}
	S && {\Sum_{s:S} p[\beta^o(s)]} && S \\
	\\
	&& S \\
	\\
	&& {S'} \\
	\\
	{S'} && {\Sum_{s':S'} p[\beta'^o(s')]} && {S'}
	\arrow["f"', from=1-1, to=7-1]
	\arrow["f"', from=3-3, to=5-3]
	\arrow[from=1-3, to=3-3]
	\arrow[from=7-3, to=5-3]
	\arrow["{\beta^u}"', from=1-3, to=1-1]
	\arrow["{\beta'^u}", from=7-3, to=7-1]
	\arrow[Rightarrow, no head, from=3-3, to=1-5]
	\arrow[Rightarrow, no head, from=5-3, to=7-5]
	\arrow["{{\beta^o}^\ast \sigma}"', from=1-5, to=1-3]
	\arrow["{{\beta'^o}^\ast \sigma}", from=7-5, to=7-3]
    \end{tikzcd}\]
    Forgetting the vertical projections out of the pullbacks gives:
    \[\begin{tikzcd}
	S && {\Sum_{s:S} p[\beta^o(s)]} && S \\
	\\
	&&&& S \\
	\\
	&&&& {S'} \\
	\\
	{S'} && {\Sum_{s':S'} p[\beta'^o(s')]} && {S'}
	\arrow["f"', from=1-1, to=7-1]
	\arrow["{\beta^u}"', from=1-3, to=1-1]
	\arrow["{\beta'^u}", from=7-3, to=7-1]
	\arrow["{{\beta^o}^\ast \sigma}"', from=1-5, to=1-3]
	\arrow["{{\beta'^o}^\ast \sigma}", from=7-5, to=7-3]
	\arrow[Rightarrow, no head, from=5-5, to=7-5]
	\arrow["f", from=3-5, to=5-5]
	\arrow[Rightarrow, no head, from=3-5, to=1-5]
    \end{tikzcd}\]
    Finally, by collapsing the identity maps and reflecting the diagram
    horizontally, we obtain
    \[\begin{tikzcd}
	S && {\Sum_{s:S} p[\beta^o(s)]} && S \\
	\\
	{S'} && {\Sum_{s':S'} p[\beta'^o(s')]} && {S'}
	\arrow["f"', from=1-1, to=3-1]
	\arrow["f", from=1-5, to=3-5]
	\arrow["{{\beta^o}^\ast \sigma}", from=1-1, to=1-3]
	\arrow["{{\beta'^o}^\ast \sigma}"', from=3-1, to=3-3]
	\arrow["{\beta^u}", from=1-3, to=1-5]
	\arrow["{\beta'^u}"', from=3-3, to=3-5]
    \end{tikzcd}\]
    which we recognize from Proposition \ref{prop:poly-dyn} as the defining
    characteristic of a morphism in \(\DynT{N}(p)\). Finally, we note that each
    of these steps is bijective, and so we have the desired bijection of
    hom-sets.
  \end{proof}
\end{prop}

\begin{question}
  Some questions:
  \begin{enumerate}%
  \item Is there a correspondence between \(\DynT{T}(p)\) and
    \(\Cat{Cat}^\#(\Fun{Cofree}_p, y^\Tt)\)?
  \item Are either of \(\DynT{T}(p)\) or \(\Cat{Cat}^\#(\Fun{Cofree}_p, y^\Tt)\)
    a topos, perhaps similarly to \(p\mdash\Cat{Coalg}\)?
  \item How do these topoi relate to behaviour topoi?
  \item How do the internal languages of these topoi relate to coalgebraic logic?
  \end{enumerate}
\end{question}

\subsection{\(pT\)-coalgebras and open Markov processes}
\label{sec:org8d59201}

In the preceding sections, we noted connections between (deterministic)
discrete-time dynamical systems over a polynomial interface \(p\) and
\(p\)-coalgebras, with Proposition \ref{prop:dyn-n-coalg} showing their
equivalence, as well as connections between random dynamical systems and Markov
chains and Markov processes (\red{REF}).  In this section, we connect these
connections by generalizing the notion of coalgebra to systems evolving in
arbitrary time.

Recall therefore that Markov chains are coalgebras for a probability monad \(\Pa
: \cat{E} \to \cat{E}\), and that open Markov chains over a polynomial interface
\(p\) are coalgebras for the composite functor \(p \Pa\), where here we
interpret \(p\) as an endofunctor on \(\cat{E}\). Recall also that a Markov
process in general time is given by a time-indexed family of Markov
kernels---\emph{i.e.}, morphisms in \(\Kl(\Pa)\)---satisfying a familiar flow
condition.  We can therefore use the recipes above to define a notion of
\(pT\)-coalgebra in general time for arbitrary monads \(T\).  Instantiating this
notion with \(T = \Pa\) will then give us a notion of open Markov process over a
polynomial interface with general time, and we can then extend the results above
to exhibit these as random dynamical systems.

\begin{defn} \label{def:pT-coalg}
  Let \(T : \cat{E} \to \cat{E}\) be a monad on the category \(\cat{E}\), and
  let \(p : \Poly{E}\) be a polynomial in \(\cat{E}\). Let \((\Tt, +, 0)\) be a
  monoid in \(\cat{E}\), representing time. Then a \(pT\)\textbf{-coalgebra with
    time} \(\Tt\) consists in a triple \(\vartheta := (S, \vartheta^o,
  \vartheta^u)\) of a \textbf{state space} \(S : \cat{E}\) and two morphisms
  \(\vartheta^o : \Tt \times S \to p(1)\) and \(\vartheta^u : \sum_{t:\Tt}
  \sum_{s:\mathbb{S}} p[\vartheta^o(t,s)] \to TS\), such that, for any section
  \(\sigma : p(1) \to \sum_{i:p(1)} p[i]\) of \(p\), the maps \(\vartheta^\sigma
  : \Tt \times S \to TS\) given by
  \[
  \Sum_{t:\Tt} S \xto{\vartheta^o(-)^\ast \sigma} \Sum_{t:\Tt} \Sum_{s:S} p[\vartheta^o(-, s)] \xto{\vartheta^u} TS
  \]
  constitute an object in the functor category \(\Cat{Cat}\big(\deloop{\Tt},
  \Kl(T)\big)\), where \(\deloop{\Tt}\) is the delooping of \(\Tt\) and
  \(\Kl(T)\) is the Kleisli category of \(T\). We call \(\vartheta^\sigma\) the
  \textbf{closure} of \(\vartheta\) by \(\sigma\).
\end{defn}

We now recall the development after Definition \ref{def:poly-dyn} of the indexed
category \(\DynT{T}\), in order to define the analogous indexed category
\(\TCoalgT{\Tt}\).

\begin{prop} \label{prop:pT-coalg}
  \(pT\) coalgebras with time \(\Tt\) form a category, denoted
  \(\PTCoalgT{\Tt}\). Its morphisms are defined as follows.  Let \(\vartheta :=
  (X, \vartheta^o, \vartheta^u)\) and \(\psi := (Y, \psi^o, \psi^u)\) be two
  \(pT\)-coalgebras. A morphism \(f : \vartheta \to \psi\) consists in a
  morphism \(f : X \to Y\) such that, for any time \(t : \Tt\) and global
  section \(\sigma : p(1) \to \Sum_{i:p(1)} p[i]\) of \(p\), the following
  naturality squares commute:
  \[\begin{tikzcd}
	X & {\Sum_{x:X} p[\vartheta^o(t, x)]} & TX \\
	Y & {\Sum_{y:Y} p[\psi^o(t, y)]} & TY
	\arrow["{\vartheta^o(t)^\ast \sigma}", from=1-1, to=1-2]
	\arrow["{\vartheta^u(t)}", from=1-2, to=1-3]
	\arrow["f"', from=1-1, to=2-1]
	\arrow["Tf", from=1-3, to=2-3]
	\arrow["{\psi^o(t)^\ast \sigma}"', from=2-1, to=2-2]
	\arrow["{\psi^u(t)}"', from=2-2, to=2-3]
  \end{tikzcd}\]
  The identity morphism \(\id_\vartheta\) on the \(pT\)-coalgebra \(\vartheta\)
  is given by the identity morphism \(\id_X\) on its state space
  \(X\). Composition of morphisms of \(pT\)-coalgebras is given by composition
  of the morphisms of the state spaces.
  \begin{proof}
    As for Proposition \ref{prop:poly-dyn}, the proof follows immediately by
    pasting.
  \end{proof}
\end{prop}

\begin{prop} \label{prop:pT-coalg-idx}
  \(\PTCoalgT{\Tt}\) extends to a polynomially-indexed category, \(\TCoalgT{\Tt}
  : \Poly{E} \to \Cat{Cat}\). Suppose \(\varphi : p \to q\) is a morphism of
  polynomials.  We define a corresponding functor \(\pTCoalgT{\varphi}{T}{\Tt} :
  \pTCoalgT{p}{T}{\Tt} \to \pTCoalgT{q}{T}{\Tt}\) as follows.  Suppose \((X,
  \vartheta^o, \vartheta^u) : \pTCoalgT{p}{T}{\Tt}\) is an object
  (\(pT\)-coalgebra) in \(\pTCoalgT{p}{T}{\Tt}\). Then
  \(\pTCoalgT{\varphi}{T}{\Tt}(X, \vartheta^o, \vartheta^u)\) is defined as the
  triple \((X, \varphi_1 \circ \vartheta^o, \vartheta^u \circ {\vartheta^o}^\ast
  \varphi^\#) : \pTCoalgT{q}{T}{\Tt}\), where the two maps are explicitly the
  following composites:
  \begin{gather*}
    \Tt \times X \xto{\vartheta^o} p(1) \xto{\varphi_1} q(1) \, ,
    \qquad
    \Sum_{t:\Tt} \Sum_{x:X} q[\varphi_1 \circ \vartheta^o(t, x)] \xto{{\vartheta^o}^\ast \varphi^\#} \Sum_{t:\Tt} \Sum_{x:X} p[\vartheta^o(t, x)] \xto{\vartheta^u} TX \, .
  \end{gather*}
  On morphisms, \(\pTCoalgT{\varphi}{T}{\Tt}(f) : \pTCoalgT{\varphi}{T}{\Tt}(X,
  \vartheta^o, \vartheta^u) \to \pTCoalgT{\varphi}{T}{\Tt}(Y, \psi^o, \psi^u)\)
  is given by the same underlying map \(f : X \to Y\) of state spaces.
  \begin{proof}
    The proof is directly analogous to that of Proposition
    \ref{prop:poly-dyn-idx}.
  \end{proof}
\end{prop}

\begin{prop} \label{prop:p-id-coalg}
  \(\PTCoalgT{\Tt}\) is equivalent to \(\DynT{T}(p)\) when \(T = \id_{\cat{E}}\).
  \begin{proof}
    This is easy to see by noting that \(\id_{\cat{E}} X = X\) for all objects
    \(X : \cat{E}\).
  \end{proof}
\end{prop}

\begin{cor}
  \(\pTCoalgT{p}{\,\id_{\cat{E}}}{\nn}\) is equivalent to \(p\mdash\Cat{Coalg}\).
  \begin{proof}
    This follows directly from Propositions \ref{prop:dyn-n-coalg} and
    \ref{prop:p-id-coalg}.
  \end{proof}
\end{cor}

\begin{rmk}
  Note that \(pT\)-coalgebras as defined above really are coalgebras \(X \to pT
  X\) in the traditional sense when \(\Tt = \nn\). We can see this by observing
  that an analogue of Proposition \ref{prop:transition-map} holds for
  \(\PTCoalgT{\nn}\), so that, following our definition, a \(pT\)-coalgebra
  \(\vartheta\) with state space \(X\) is determined by two morphisms
  \(\vartheta^o : X \to p(1)\) and \(\vartheta^u : \sum_{x:X} p[\vartheta^o(x)]
  \to TX\). Then, note that a `classical' \(pT\)-coalgebra \(\vartheta' : X \to
  pT X\) is equivalently a morphism \(\vartheta : X \to \sum_{i:p(1)}
  X^{p[i]}\), by the definition of the polynomial functor \(p\). But, by the
  universal property of the dependent sum (and as in the proof of Proposition
  \ref{prop:dyn-n-coalg}), such a morphism corresponds bijectively to such a
  pair of maps \((\vartheta^o, \vartheta^u)\) as determines our earlier
  \(pT\)-coalgebra \(\vartheta\)! And as in Proposition \ref{prop:dyn-n-coalg},
  our definition of \(pT\)-coalgebra morphism corresponds to the classical
  notion of coalgebra homomorphism, so that our category \(\PTCoalgT{\nn}\) is
  equivalent to the classical category \(\Cat{Coalg}(pT)\) of \(pT\)-coalgebras
  and coalgebra homomorphisms. This justifies our thinking of the categories
  \(\PTCoalgT{\Tt}\) as generalized categories of coalgebras.
\end{rmk}

\begin{prop}
  \(\pTCoalgT{y}{T}{\Tt}\) is equivalent to \(\Cat{Cat}\big(\deloop{\Tt},
  \Kl(T)\big)\).
  \begin{proof}
    Analogous to the proof of Proposition \ref{prop:closed-sys-in-dyn-y} for the
    equivalence \(\DynT{T}(y) \cong \Cat{Cat}(\deloop{\Tt}, \cat{E})\).
  \end{proof}
\end{prop}

\begin{rmk}[Closed Markov chains and Markov processes]
  A closed \textit{Markov chain} is given by a map \(X \to \Pa X\), where \(\Pa
  : \cat{E} \to \cat{E}\) is a probability monad on \(\cat{E}\); this is
  equivalently a \(y\Pa\)-coalgebra with time \(\nn\), and an object in
  \(\Cat{Cat}\big(\deloop{\nn}, \Kl(\Pa)\big)\). With more general time \(\Tt\),
  one obtains closed \textit{Markov processes}: objects in
  \(\Cat{Cat}\big(\deloop{\Tt}, \Kl(\Pa)\big)\). More explicitly, a closed
  Markov process is a time-indexed family of Markov kernels; that is, a morphism
  \(\vartheta : \Tt \times X \to \Pa X\) such that, for all times \(s,t : \Tt\),
  \(\vartheta_{s+t} = \vartheta_s \klcirc \vartheta_t\) as a morphism in
  \(\Kl(\Pa)\). Note that composition \(\klcirc\) in \(\Kl(\Pa)\) is given by
  the Chapman-Kolmogorov equation, so this means that
  \[
  \vartheta_{s+t}(y|x) = \int_{x':X} \vartheta_s(y|x') \, \vartheta_t(\d x'|x) \, .
  \]
\end{rmk}

\emph{Open} Markov processes over a polynomial interface \(p\) are therefore the
objects of the category \(\pTCoalgT{p}{\Pa}{\Tt}\) for a given time monoid
\(\Tt\); the generalized flow condition corresponds to the satisfaction of an
analogous Chapman-Kolmogorov equation by the closures of the systems by any
section of \(p\).

We now translate the development of bundle and nested open dynamical systems
(after Definition \ref{defn:bdyn-pbth}) to the coalgebraic setting. Our
principal aim is to define a notion of ``nested Markov process'' and
corresponding generalizations for arbitrary monads \(T\).

\begin{defn} \label{defn:bun-coalg}
  Let \(p, b : \Poly{E}\) be polynomials in \(\cat{E}\), and let \(\theta :=
  (\theta(\ast), \theta^o, \theta^u) : \pTCoalgT{b}{T}{\Tt}\) be a
  \(bT\)-coalgebra.  A \(pT\)\textbf{-coalgebra bundle} over \(\theta\) is a
  pair \((\pi_{\vartheta\theta}, \vartheta)\) where \(\vartheta :=
  (\vartheta(\ast), \vartheta^o, \vartheta^u) : \pTCoalgT{p}{T}{\Tt}\) is a
  \(pT\)-coalgebra and \(\pi_{\vartheta\theta} : \vartheta(\ast) \to
  \theta(\ast)\) is a bundle in \(\cat{E}\), such that, for all time \(t : \Tt\)
  and sections \(\sigma\) of \(p\) and \(\varsigma\) of \(b\), the following
  diagrams commute, thereby inducing a bundle of closed dynamical systems
  \(\pi^{\sigma\varsigma}_{\vartheta\theta} : \vartheta^\sigma \to
  \theta^\varsigma\) in \(\Cat{Cat}\big(\deloop{\Tt}, \Kl(T)\big)\):
  \[\begin{tikzcd}
	\vartheta(\ast) && {\Sum_{w:\vartheta(\ast)} p[\vartheta^o(t, w)]} && {T \vartheta(\ast)} \\
	\\
	\theta(\ast) && {\Sum_{x:\theta(\ast)} b[\theta^o(t,x)]} && {T \theta(\ast)}
	\arrow["{\pi_{\vartheta\theta}}", from=1-1, to=3-1]
	\arrow["{T\pi_{\vartheta\theta}}", from=1-5, to=3-5]
	\arrow["{\vartheta^o(t)^\ast \sigma}", from=1-1, to=1-3]
	\arrow["{\vartheta^u(t)}", from=1-3, to=1-5]
	\arrow["{\theta^o(t)^\ast \varsigma}", from=3-1, to=3-3]
	\arrow["{\theta^u(t)}", from=3-3, to=3-5]
  \end{tikzcd}\]
\end{defn}

\begin{prop} \label{prop:bun-coalg-cat}
  Let \(p, b : \Poly{E}\) be polynomials in \(\cat{E}\), and let \(\theta :=
  (\theta(\ast), \theta^o, \theta^u) : \pTCoalgT{b}{T}{\Tt}\) be a
  \(bT\)-coalgebra.  \(pT\)-coalgebra bundles over \(\theta\) form the objects
  of a category \(\pTCoalgT{(p,b)}{T}{\Tt}/\theta\). Morphisms \(f :
  (\pi_{\vartheta\theta}, \vartheta) \to (\pi_{\varrho\theta}, \varrho)\) are
  maps \(f : \vartheta(\ast) \to \varrho(\ast)\) in \(\cat{E}\) making the
  following diagram commute for all times \(t : \Tt\) and sections \(\sigma\) of
  \(p\) and \(\varsigma\) of \(b\):
  \[\begin{tikzcd}
	{\vartheta(\ast)} &&&& {\Sum_{w:\vartheta(\ast)} p[\vartheta^o(t, w)]} &&&& {T\vartheta(\ast)} \\
	\\
	&& {\theta(\ast)} && {\Sum_{x:\theta(\ast)} b[\theta^o(t,x)]} && {T\theta(\ast)} \\
	\\
	{\varrho(\ast)} &&&& {\Sum_{y:\varrho(\ast)} p[\varrho^o(t, y)]} &&&& {T\varrho(\ast)}
	\arrow["{\vartheta^o(t)^\ast \sigma}", from=1-1, to=1-5]
	\arrow["{\vartheta^u(t)}", from=1-5, to=1-9]
	\arrow["{\theta^o(t)^\ast \varsigma}", from=3-3, to=3-5]
	\arrow["{\theta^u(t)}", from=3-5, to=3-7]
	\arrow["{\pi_{\vartheta\theta}}"', from=1-1, to=3-3]
	\arrow["{\pi_{\varrho\theta}}", from=5-1, to=3-3]
	\arrow["{T\pi_{\vartheta\theta}}", from=1-9, to=3-7]
	\arrow["{T\pi_{\varrho\theta}}"', from=5-9, to=3-7]
	\arrow["{\varrho^o(t)^\ast \sigma}"', from=5-1, to=5-5]
	\arrow["{\varrho^u(t)}"', from=5-5, to=5-9]
	\arrow["{f}"', from=1-1, to=5-1]
	\arrow["{Tf}", from=1-9, to=5-9]
  \end{tikzcd}\]
  That is, \(f\) is a map on the state spaces that induces a morphism
  \((\pi_{\vartheta\theta}, \vartheta^\sigma) \to (\pi_{\varrho\theta},
  \varrho^\sigma)\) in \(\Cat{Cat}(\deloop{\Tt}, \Kl(T))/\theta^\varsigma\) of
  bundles of the closures. Identity morphisms are the corresponding identity
  maps, and composition is by pasting.
\end{prop}

\begin{prop} \label{prop:bun-coalg-bth-idx}
  Varying the polynomials \(p\) in \(\pTCoalgT{(p,b)}{T}{\Tt}/\theta\) induces a
  polynomially indexed category \(\pTCoalgT{(-,b)}{T}{\Tt}/\theta : \Poly{E} \to
  \Cat{Cat}\). On polynomials \(p\), it returns the categories
  \(\pTCoalgT{(p,b)}{T}{\Tt}/\theta\) of Proposition
  \ref{prop:bun-coalg-cat}. On morphisms \(\varphi : p \to q\) of polynomials,
  define the functors \(\pTCoalgT{(\varphi,b)}{T}{\Tt}/\theta :
  \pTCoalgT{(p,b)}{T}{\Tt}/\theta \to \pTCoalgT{(q,b)}{T}{\Tt}/\theta\) as in
  Proposition \ref{prop:pT-coalg-idx}. That is, suppose
  \((\pi_{\vartheta\theta}, \vartheta) : \pTCoalgT{(p,b)}{T}{\Tt}/\theta\) is an
  object in \(\pTCoalgT{(p,b)}{T}{\Tt}/\theta\), where \(\vartheta :=
  (\vartheta(\ast), \vartheta^o, \vartheta^u)\). Then its image
  \(\left(\pTCoalgT{(\varphi,b)}{T}{\Tt}/\theta\right)(\pi_{\vartheta\theta},
  \vartheta)\) is defined as the pair \((\pi_{\vartheta\theta},
  \varphi\vartheta)\), where \(\varphi\vartheta := (\vartheta(\ast), \phi_1
  \circ \vartheta^o, \vartheta^u \circ {\vartheta^o}^\ast \varphi^\#)\). On
  morphisms \(f : (\pi_{\vartheta\theta}, \vartheta) \to (\pi_{\varrho\theta},
  \varrho)\), \(\left(\pTCoalgT{(\varphi,b)}{T}{\Tt}/\theta\right)(f)\) is again
  given by the same underlying map \(f : \vartheta(\ast) \to \varrho(\ast)\) of
  state spaces.
  \begin{proof}
    Analogous to the proof of Proposition \ref{prop:bdyn-bth-idx}.
  \end{proof}
\end{prop}

\begin{prop} \label{prop:bun-coalg-b-idx}
  Letting the base system \(\theta\) also vary induces a doubly-indexed category
  \(\pTCoalgT{({-},b)}{T}{\Tt}/({=}) : \Poly{E} \times \pTCoalgT{b}{T}{\Tt} \to
  \Cat{Cat}\). Given a polynomial \(p : \Poly{E}\) and morphism \(\phi : \theta
  \to \rho\) in \(\pTCoalgT{b}{T}{\Tt}\), the functor
  \(\pTCoalgT{(p,b)}{T}{\Tt}/\phi : \pTCoalgT{(p,b)}{T}{\Tt}/\theta \to
  \pTCoalgT{(p,b)}{T}{\Tt}/\rho\) is defined by post-composition, as in
  Propositions \ref{prop:poly-rdyn-idx} and \ref{prop:bdyn-b-idx}. More
  explicitly, such a morphism \(\phi\) corresponds to a map \(\phi :
  \theta(\ast) \to \rho(\ast)\) of state spaces in \(\cat{E}\). Given an object
  \((\pi_{\vartheta\theta}, \vartheta)\) of \(\pTCoalgT{(p,b)}{T}{\Tt}/\theta\),
  we define \(\left(\pTCoalgT{(p,b)}{T}{\Tt}/\phi\right)(\pi_{\vartheta\theta},
  \vartheta) := (\phi \circ \pi_{\vartheta\theta}, \vartheta)\). Given a
  morphism \(f : (\pi_{\vartheta\theta}, \vartheta) \to (\pi_{\varrho\theta},
  \varrho)\) in \(\pTCoalgT{(p,b)}{T}{\Tt}/\theta\), its image
  \(\left(\pTCoalgT{(p,b)}{T}{\Tt}/\phi\right)(f) :
  (\phi\circ\pi_{\vartheta\theta}, \vartheta) \to (\phi\circ\pi_{\varrho\theta},
  \varrho)\) is given by the same underlying map \(f : \vartheta(\ast) \to
  \varrho(\ast)\) of state spaces.
  \begin{proof}
    Analogous to the proof of Proposition \ref{prop:bdyn-b-idx}.
  \end{proof}
\end{prop}

\begin{prop} \label{prop:bun-coalg-b-fib}
  There is an indexed opfibration \(\int\pTCoalgT{(-,b)}{T}{\Tt} : \Poly{E} \to
  \Cat{Fib}\left(\pTCoalgT{b}{T}{\Tt}\right)\) generated from
  \(\pTCoalgT{({-},b)}{T}{\Tt}/({=})\) by the Grothendieck construction.

  Explicitly, given a polynomial \(p : \Poly{E}\), the objects of
  \(\int\pTCoalgT{(p,b)}{T}{\Tt}\) are triples \((\pi_{\vartheta\theta},
  \vartheta, \theta)\), where \(\theta : \pTCoalgT{b}{T}{\Tt}\) is a
  \(bT\)-coalgebra and \((\pi_{\vartheta\theta}, \vartheta) :
  \pTCoalgT{(p,b)}{T}{\Tt}/\theta\) is a \(pT\)-coalgebra bundle over
  \(\theta\). Morphisms \(f : (\pi_{\vartheta\theta}, \vartheta, \theta) \to
  (\pi_{\varrho\rho}, \varrho, \rho)\) are pairs \((f_p, f_b)\) of a morphism
  \(f_b : \theta \to \rho\) in \(\pTCoalgT{b}{T}{\Tt}\) and a morphism \(f_p :
  (f_b \circ \pi_{\vartheta\theta}, \vartheta) \to (\pi_{\varrho\rho},
  \varrho)\) in \(\pTCoalgT{(p,b)}{T}{\Tt}/\theta\) making the following diagram
  commute for all sections \(\sigma\) of \(p\) and \(\varsigma\) of \(b\):
  \[\begin{tikzcd}
	{\vartheta(\ast)} && {\Sum_{w:\vartheta(\ast)} p[\vartheta^o(t, w)]} && {T\vartheta(\ast)} \\
	\\
	& {\theta(\ast)} && {\Sum_{x:\theta(\ast)} b[\theta^o(t, x)]} && {T\theta(\ast)} \\
	\\
	{\varrho(\ast)} && {\Sum_{y:\varrho(\ast)} p[\varrho^o(t,y)]} && {T\varrho(\ast)} \\
	\\
	& {\rho(\ast)} && {\Sum_{z:\rho(\ast)} b[\rho^o(t, z)]} && {T\rho(\ast)}
	\arrow["{\pi_{\vartheta\theta}}", from=1-1, to=3-2]
	\arrow["{T\pi_{\vartheta\theta}}", from=1-5, to=3-6]
	\arrow["{\vartheta^o(t)^\ast \sigma}", from=1-1, to=1-3]
	\arrow["{\theta^o(t)^\ast \varsigma}", from=3-2, to=3-4]
	\arrow["{\vartheta^u(t)}", from=1-3, to=1-5]
	\arrow["{\theta^u(t)}"{pos=0.3}, from=3-4, to=3-6]
	\arrow["{\varrho^o(t)^\ast \sigma}"{pos=0.7}, from=5-1, to=5-3]
	\arrow["{\varrho^u(t)}", from=5-3, to=5-5]
	\arrow["{\rho^o(t)^\ast \varsigma}", from=7-2, to=7-4]
	\arrow["{\rho^u(t)}", from=7-4, to=7-6]
	\arrow["{\pi_{\varrho\rho}}"', from=5-1, to=7-2]
	\arrow["{T\pi_{\varrho\rho}}"', from=5-5, to=7-6]
	\arrow["{f_p}"'{pos=0.4}, curve={height=6pt}, from=1-1, to=5-1]
	\arrow["{f_b}"'{pos=0.3}, curve={height=6pt}, from=3-2, to=7-2]
	\arrow["{T f_p}"{pos=0.7}, curve={height=-6pt}, from=1-5, to=5-5]
	\arrow["{T f_b}"{pos=0.6}, curve={height=-6pt}, from=3-6, to=7-6]
  \end{tikzcd}\]
  Identity morphisms are the pairs of the corresponding identities, and
  composition is again by pasting.
  \begin{proof}
    Compare Propositions \ref{prop:bdyn-b-fib} and \ref{prop:poly-rdyn-fib}.
  \end{proof}
\end{prop}

\begin{prop} \label{prop:bun-coalg-idx}
  Varying the base polynomial \(b\) extends \(\pTCoalgT{({-},b)}{T}{\Tt}/({=})\)
  to a triply indexed category, \(\pTCoalgT{({-},{=})}{T}{\Tt}/({\equiv}) :
  \Poly{E} \times \sum_{b:\Poly{E}} \pTCoalgT{b}{T}{\Tt} \to \Cat{Cat}\) and
  \(\int\pTCoalgT{({-},b)}{T}{\Tt}\) to a doubly indexed fibration,
  \(\int\pTCoalgT{({-},{=})}{T}{\Tt} : \Poly{E} \to \prod_{b:\Poly{E}}
  \Cat{Fib}\left(\pTCoalgT{b}{T}{\Tt}\right)\); these are equivalent by the
  Grothendieck construction.

  Let \(p, b, c : \Poly{E}\) be polynomials, and let \(\chi : b \to c\) be a
  morphism accordingly. Let \(\theta : \pTCoalgT{b}{T}{\Tt}\) range over
  \(bT\)-coalgebras, and let \(\chi_\ast := \pTCoalgT{\chi}{T}{\Tt}\). We define
  the (dependent) functor \(\pTCoalgT{(p,\chi)}{T}{\Tt}/\theta :
  \pTCoalgT{(p,b)}{T}{\Tt}/\theta \allowbreak\to
  \pTCoalgT{(p,c)}{T}{\Tt}/\chi_\ast\theta\) as follows.  This functor is
  equivalent to its image under the Grothendieck construction,
  \(\int\pTCoalgT{(p,\chi)}{T}{\Tt} : \int\pTCoalgT{(p,b)}{T}{\Tt} \to
  \int\pTCoalgT{(p,c)}{T}{\Tt}\), which is easier to describe. Therefore, let
  \((\pi_{\vartheta\theta}, \vartheta, \theta)\) be an object of
  \(\int\pTCoalgT{(p,b)}{T}{\Tt}\). Its image under
  \(\int\pTCoalgT{(p,\chi)}{T}{\Tt}\) is the object \((\pi_{\vartheta\theta},
  \vartheta, \chi_\ast \theta) : \int\pTCoalgT{(p,c)}{T}{\Tt}\). Given a
  morphism \((f_p, f_b) : (\pi_{\vartheta\theta}, \vartheta, \theta) \to
  (\pi_{\varrho\rho}, \varrho, \rho)\) in \(\int\pTCoalgT{(p,b)}{T}{\Tt}\), its
  image \(\left(\int\pTCoalgT{(p,\chi)}{T}{\Tt}\right)(f) :
  (\pi_{\vartheta\theta}, \vartheta, \chi_\ast \theta) \to (\pi_{\varrho\rho},
  \varrho, \chi_\ast \rho)\) is given by the same maps \(f_p : \vartheta(\ast)
  \to \varrho(\ast)\) and \(f_b : \theta(\ast) \to \rho(\ast)\) of state spaces.
  \begin{proof}
    Compare Proposition \ref{prop:bdyn-idx}.
  \end{proof}
\end{prop}

And, when \(p\) is nested over \(b\) we of course have a correspondingly nested
notion of coalgebra:

\begin{prop} \label{prop:ncoalg-idx}
  There is an indexed category \(\pTCoalgT{({-},{=})}{T}{\Tt}/({\equiv}) :
  \sum_{b:\Poly{E}(y)} \Poly{E}(b) \times \pTCoalgT{b}{T}{\Tt} \to \Cat{Cat}\)
  of \textbf{generalized coalgebras over nested polynomials}, and a
  corresponding dependently-indexed fibration \(\int\pTCoalgT{({-},{=})}{T}{\Tt}
  : \prod_{b:\Poly{E}(y)} \Poly{E}(b) \to
  \Cat{Fib}\big(\pTCoalgT{b}{T}{\Tt}\big)\). We have purposefully overloaded the
  names of these categories, since they are defined as the categories of
  coalgebra bundles, only with an extra condition of compatibility with the
  nesting. That is, suppose that \(p\) is a polynomial nested over \(b\). Then
  there are morphisms \(m : \sum_i p[i] \to \sum_k b[k]\) and \(n : p(1) \to
  b(1)\) in \(\cat{E}\), such that \(n \circ p = b \circ m\). A generalized
  \(T\)-coalgebra over \((m, n) : p \to b\), \textit{i.e.} an object of
  \(\int\pTCoalgT{(m,n)}{T}{\Tt}\), is a triple \((\pi_{\vartheta\theta},
  \vartheta, \theta)\) where \(\theta : \pTCoalgT{b}{T}{\Tt}\) is a
  \(bT\)-coalgebra, \(\vartheta : \pTCoalgT{p}{T}{\Tt}\) is a \(pT\)-coalgebra,
  and \(\pi_{\vartheta\theta} : \vartheta(\ast) \to \theta(\ast)\) is a map
  between the state spaces which lifts uniquely to make the front face of the
  following cube, and therefore the whole cube, commute:
  \[\begin{tikzcd}
	&& {\Sum_{i:p(1)} p[i]} &&& {\Sum_{k:b(1)} b[k]} \\
	\\
	{\Sum_{w:\vartheta(\ast)} p[\vartheta(t, w)]} &&& {\Sum_{x:\theta(\ast)} b[\theta(t,x)]} \\
	& {} & {p(1)} &&& {b(1)} \\
	\\
	{\vartheta(\ast)} &&& {\theta(\ast)}
	\arrow["{\pi_{\vartheta\theta}}", from=6-1, to=6-4]
	\arrow[dashed, from=3-1, to=3-4]
	\arrow[from=3-1, to=6-1]
	\arrow[from=3-4, to=6-4]
	\arrow[from=3-1, to=1-3]
	\arrow[from=3-4, to=1-6]
	\arrow["m", from=1-3, to=1-6]
	\arrow[from=1-3, to=4-3]
	\arrow[from=1-6, to=4-6]
	\arrow["n", from=4-3, to=4-6]
	\arrow["{\vartheta^o(t)}"'{pos=0.6}, from=6-1, to=4-3]
	\arrow["{\theta^o(t)}"', from=6-4, to=4-6]
  \end{tikzcd}\]
  As in the dynamical-system case (Proposition \ref{prop:ndyn-idx}), we call the
  above condition the \textit{nesting condition}; the unlabelled edges of this
  cube are the obvious projections.  As usual, we also require a flow condition
  to be satisfied; that is, for all times \(t : \Tt\) and sections \(\sigma\) of
  \(p\) and \(\varsigma\) of \(b\), the following diagrams commute, where the
  dashed arrow below is the same dashed lift above:
  \[\begin{tikzcd}
	{\vartheta(\ast)} && {\Sum_{w:\vartheta(\ast)} p[\vartheta^o(t,w)]} && {T\vartheta(\ast)} \\
	\\
	{\theta(\ast)} && {\Sum_{x:\theta(\ast)} b[\theta^o(t,x)]} && {T\theta(\ast)}
	\arrow["{\pi_{\vartheta\theta}}", from=1-1, to=3-1]
	\arrow["{T\pi_{\vartheta\theta}}", from=1-5, to=3-5]
	\arrow["{\vartheta^o(t)^\ast \sigma}", from=1-1, to=1-3]
	\arrow["{\sigma^o(t)^\ast \varsigma}"', from=3-1, to=3-3]
	\arrow["{\vartheta^u(t)}", from=1-3, to=1-5]
	\arrow["{\theta^u(t)}"', from=3-3, to=3-5]
	\arrow[dashed, from=1-3, to=3-3]
  \end{tikzcd}\]
  Morphisms in \(\int\pTCoalgT{(m,n)}{T}{\Tt}\) are as for
  \(\int\pTCoalgT{(p,b)}{T}{\Tt}\) (Proposition \ref{prop:bun-coalg-idx}), with
  the addition of the dashed arrows on the top and bottom of the defining cube
  (in \ref{prop:bun-coalg-idx}). The actions of
  \(\pTCoalgT{({-},{=})}{T}{\Tt}/({\equiv})\) and of
  \(\int\pTCoalgT{({-},{=})}{T}{\Tt}\) on morphisms of polynomials, nested
  polynomials, and base coalgebras are defined as in the non-nested case above;
  given a morphism of polynomials, the dashed lifts are transformed by
  pullback. Because the pasting of two pullback squares is again a pullback
  square, it is easy to check that this also constitutes an indexed category.
\end{prop}

\begin{question}
  Some questions:

  \begin{enumerate}
  \item Under what conditions can we say something like,
    ``\(\pTCoalgT{p}{T}{\Tt}\) is equivalent to \(\Cat{Cat}\big(\deloop{\Tt},
    \Kl(pT)\big)\)''?

  \item What about coalgebraic logic? Distributive laws? Does that good stuff
    lift to this regime?

  \item What about mapping cospan-algebra processes to processes over polynomials?
  \end{enumerate}
\end{question}

\section{References}
\label{sec:org5aae1c2}
\printbibliography[heading=none]

\appendix
\end{document}